\documentclass{amsart}
\usepackage{amssymb,latexsym}
\usepackage{amsthm,amsmath,amscd}
\newcommand{\co}{\colon\thinspace}

\newtheorem{prop}{Proposition}
\newtheorem{lemma}[prop]{Lemma}
\newtheorem{thm}[prop]{Theorem}
\newtheorem{cor}[prop]{Corollary}
\newtheorem{sub}[prop]{Sublemma}
\newtheorem{dfn}[prop]{Definition}
\newtheorem{ass}[prop]{Assumption}
\newtheorem{notn}[prop]{Notation}
\newtheorem{conj}[prop]{Conjecture}
\theoremstyle{remark}
\newtheorem{rem}[prop]{Remark}
 
\DeclareMathOperator{\ep}{\epsilon}
\DeclareMathOperator{\dbar}{\bar{\partial}}
\numberwithin{equation}{section} \numberwithin{prop}{section}
\title{Standard surfaces and nodal curves in symplectic 4-manifolds}
\author{Michael Usher}
\email{usher@alum.mit.edu}
\begin{document}
\begin{abstract}
Continuing the program of \cite{DS} and \cite{Usher}, we introduce
refinements of the Donaldson-Smith standard surface count which
are designed to count nodal pseudoholomorphic curves and curves
with a prescribed decomposition into reducible components.  In
cases where a corresponding analogue of the Gromov-Taubes
invariant is easy to define, our invariants agree with those
analogues.  We also prove a vanishing result for some of the
invariants that count nodal curves.
\end{abstract}
\maketitle

\section{Introduction}

Let $(X,\omega)$ be a closed symplectic 4-manifold.  We assume
that $[\omega]\in H^2(X,\mathbb{Z})$; however, the main theorems
in this paper concern Gromov invariants, which are unchanged under
deformations of the symplectic form, so since any symplectic form
is deformation equivalent to an integral form there is no real
loss of generality here.  According to \cite{Don}, if $k$ is large
enough, taking a suitable pair of sections of a line bundle
$L^{\otimes k}$ where $L$ has Chern class $[\omega]$ and blowing
$X$ up at the common vanishing locus of these sections to obtain
the new manifold $X'$ gives rise to a symplectic Lefschetz
fibration $f\co X'\to \mathbb{C}P^1$ (the exceptional curves of
the blowup $\pi\co X'\to X$ appear as sections of $f$, while at
other points $x'\in X'$, $f(x')\in\mathbb{C}\cup\{\infty\}$ is the
ratio of the two chosen sections of $L^{\otimes k}$ at $\pi
(x')\in X$). In other words, $f$ is a fibration by Riemann
surfaces over the complement of a finite set of critical values in
$S^2$, while near its critical points $f$ is given in smooth local
complex coordinates by $f(z,w)=zw$.  Results of \cite{Smith2} show
that the critical points of $f$ may be assumed to lie in separate
fibers, and all fibers of $f$ may be assumed irreducible.  Once we
choose a metric on $X'$, Donaldson's construction thus presents a
suitable blowup of $X$ as a smoothly $\mathbb{C}P^1$-parametrized
family of Riemann surfaces, all but finitely of which are smooth
and all of which are irreducible with at worst one ordinary double
point. Where $\kappa_X=c_1 (T^* X)$ is the canonical class of $X$,
note that the adjunction formula gives the arithmetic genus of the
fibers as $g=1+(k^2[\omega]^2+k\kappa_X\cdot\omega)/2$.

Beginning with the work of S. Donaldson and I. Smith in \cite{DS},
some efforts have recently been made toward determining whether
such a Lefschetz fibration can shed light on any questions
concerning pseudoholomorphic curves in $X$.  More specifically,
for any natural number $r$ Donaldson and Smith construct the
\emph{relative Hilbert scheme}, which is a smooth symplectic
manifold $X_r (f)$ with a map $F\co X_r (f)\to \mathbb{C}P^1$
whose fiber over a regular value $t$ of $f$ is the symmetric
product $S^r f^{-1}(t)$.  If we choose an almost complex structure
$j$ on $X'$ with respect to which $f$ is a pseudoholomorphic map,
a $j$-holomorphic curve $C$ in $X'$ which contains no fiber
components will, by the positivity of intersections between
$j$-holomorphic curves, meet each fiber in $r:=[C]\cdot [fiber]$
points, counted with multiplicities.  In other words, $C\cap
f^{-1}(t)\in S^{r}f^{-1}(t)$, so that, letting $t$ vary, $C$ gives
rise to a section $s_C$ of $X_r (f)$. Conversely, a section $s$ of
$X_{r}(f)$ gives rise to a subset $C_s$ of $X'$ (namely the union
of all the points appearing in the divisors $s(t)$ as $t$ varies),
and from $j$ one may construct a (nongeneric and generally not
even $C^1$) almost complex structure $\mathbb{J}_j$ with the
property that $C$ is a (possibly disconnected) $j$-holomorphic
curve in $X'$ if and only if $s_C$ is a $\mathbb{J}_j$-holomorphic
section of $X_{r}(f)$.

Accordingly, it seems reasonable to study pseudoholomorphic curves
in $X'$ by studying pseudoholomorphic sections of $X_{r}(f)$.  If
$\alpha\in H^2 (X';\mathbb{Z})$, the \emph{standard surface count}
$\mathcal{DS}_f (\alpha)$ is defined in \cite{Smith} (and earlier
in \cite{DS} for $\alpha=\kappa_{X'}$) as the Gromov-Witten
invariant which counts $J$-holomorphic sections $s$ whose
corresponding sets $C_s$ are Poincar\'e dual to the class $\alpha$
and pass through a generic set of
$d(\alpha)=\frac{1}{2}(\alpha^2-\kappa_{X'}\cdot \alpha)$ points
of $X'$, where $J$ is a generic almost complex structure on $X_r
(f)$. Smith shows in \cite{Smith} that there is at most one
homotopy class $c_{\alpha}$ of sections $s$ such that $C_s$ is
Poincar\'e dual to $\alpha$, and moreover that the complex
dimension of the space of $J$-holomorphic sections in this
homotopy class is, for generic $J$, the aforementioned
$d(\alpha)$, which the reader may recognize as the same as the
expected dimension of $j$-holomorphic submanifolds of $X$
Poincar\'e dual to $\alpha$. Further, the moduli space of
$J$-holomorphic sections in the homotopy class $c_{\alpha}$ is
compact for generic $J$ if $k$ is taken large enough.  The moduli
space in the definition of $\mathcal{DS}_f$ is therefore a finite
set, and $\mathcal{DS}_f$ simply counts the members of this set
with sign according to the usual (spectral-flow-based)
prescription.

Donaldson and Smith have proven various results about
$\mathcal{DS}$, perhaps the most notable of which is the main
theorem of \cite{Smith}, which asserts that if $\alpha\in
H^{2}(X;\mathbb{Z})$, if $b^{+}(X)>b_1 (X)+1$, and if the degree
$k$ of the Lefschetz fibration is high enough, then
\begin{equation} \label{SD} \mathcal{DS}_f (\pi^* \alpha)=\pm \mathcal{DS}_f (\pi^*
(\kappa_X - \alpha)). \end{equation} Their work has led to new,
more symplectic proofs of various results in 4-dimensional
symplectic topology which had previously been accessible only by
Seiberg-Witten theory (as an example we mention the main theorem
of \cite{DS}, according to which $X$ admits a symplectic surface
Poincar\'e dual to $\kappa_X$, again assuming $b^{+}(X)>b_1
(X)+1$).  In \cite{Usher} it was shown that the invariant
$\mathcal{DS}_f$ agrees with the Gromov invariant $Gr$ which was
introduced  by C. Taubes in \cite{Taubes} and which counts
possibly-disconnected pseudoholomorphic submanifolds of $X'$ Poincar\'e dual to a
given cohomology class.  This in particular shows that
$\mathcal{DS}_f$ is independent of the choice of Lefschetz
fibration structure, and, in combination with Smith's duality
theorem (\ref{SD}) and the fact that under a blowup $\pi$ one has
$Gr(\pi^* \alpha)=Gr(\alpha)$, yields a new proof of the relation
\[ Gr(\alpha)=\pm Gr(\kappa_X-\alpha)\] if $b^{+}(X)>b_1 (X)+1$, a
result which had previously only been known as a shadow of the
charge conjugation symmetry in Seiberg-Witten theory.

The information contained in the Gromov invariants comprises only
a part of the data that might be extracted from pseudoholomorphic
curves in $X$.  The present paper aims to show that many of these
additional data can also be captured by Donaldson-Smith-type
invariants.  For instance, $Gr(\alpha)$ counts all of the curves
with any decomposition into connected components whose homology
classes add up (counted with multiplicities) to $\alpha$. It is
natural to wish to keep track of the decompositions of our curves
into reducible components; accordingly we make the following:

\begin{dfn} \label{invts} Let $\alpha\in H^{2}(X;\mathbb{Z})$.  Let \[
\alpha=\beta_1+\cdots +\beta_m+c_1 \tau_1+\cdots +c_n \tau_n \] be
a decomposition of $\alpha$ into distinct summands, where none of
the $\beta_i$ satisfies $\beta_{i}^{2}=\kappa_X \cdot\beta_i =0$,
while the $\tau_i$ are distinct classes which are primitive in the
lattice $H^{2}(X;\mathbb{Z})$ and all satisfy $\tau_{i}^2=\kappa_X
\cdot \tau_i =0$.  Then \[
Gr(\alpha;\beta_1,\ldots,\beta_m,c_1\tau_1,\cdots,c_n\tau_n)
\] is the invariant counting ordered $(m+n)$-tuples $(C_1,\ldots,C_{m+n})$ of
transversely intersecting smooth pseudoholomorphic curves in $X$,
where
\begin{itemize} \item[(i)] for $1\leq i\leq m$, $C_i$ is a connected curve Poincar\'e dual to $\beta_i$
which passes through some prescribed generic set of $d(\beta_i)$
points; \item[(ii)] for $m+1\leq k\leq m+n$, $C_k$ is a union of
connected curves Poincar\'e dual to classes
$l_{k,1}\tau_k,\cdots,l_{k,p}\tau_k$ decorated with positive
integer multiplicities $m_{k,q}$ with the property that \[\sum_q
m_{k,q}l_{k,q}=c_k.\] \end{itemize} The weight of each component
of each such curve is to be determined according to the
prescription given in the definition of the Gromov invariant in
\cite{Taubes} (in particular, the components $C_{k,q}$ in class
$l_{k,q}\tau_k$ are given the weight $r(C_{k,q},m_{k,q})$
specified in Section 3 of \cite{Taubes}), and the contribution of
the entire curve is the product of the weights of its components.
\end{dfn}

The objects counted by $Gr(\alpha;\alpha_1,\ldots,\alpha_n)$  will
then be reducible curves with smooth irreducible components and a
total of $\sum \alpha_i\cdot \alpha_j$ nodes arising from
intersections between these components. $Gr(\alpha)$ is the sum
over all decompositions of $\alpha$ into classes which are
pairwise orthogonal under the cup product of the
\[ \frac{d(\alpha)!}{\prod (d(\alpha_i)!)}Gr(\alpha;\alpha_1,\ldots,\alpha_n);\] in turn, one has
\[
Gr(\alpha;\alpha_1,\ldots,\alpha_n)=\prod_{i=1}^{n}Gr(\alpha_i;\alpha_i).\]

In Section \ref{red}, given a symplectic Lefschetz fibration $f\co
X\to S^2$ with sufficiently large fibers, by counting sections of
a relative Hilbert scheme we construct a corresponding invariant
$\widetilde{\mathcal{DS}}_f(\alpha;\alpha_1,\ldots,\alpha_n)$
provided that none of the $\alpha_i$ can be written as $m\beta$
where $m>1$ and $\beta$ is Poincar\'e dual to either a symplectic
square-zero torus or a symplectic $(-1)$-sphere.  Further:

\begin{thm}$\frac{(\sum
d(\alpha_i))!}{\prod
(d(\alpha_i)!)}Gr(\alpha;\alpha_1,\ldots,\alpha_n)=\widetilde{\mathcal{DS}_f}(\alpha;\alpha_1,\ldots,\alpha_n)$
provided that the degree of the fibration is large enough that
$\langle [\omega_{X}],fiber\rangle>[\omega_{X}]\cdot
\alpha$.\end{thm}

The sections $s$ counted by
$\widetilde{\mathcal{DS}}_f(\alpha;\alpha_1,\ldots,\alpha_n)$
correspond tautologically to curves $C_s=\cup C^{i}_{s}$ in $X$
with each $C^{i}_{s}$ Poincar\'e dual to $\alpha_i$.  The
$C_{s}^{i}$ will be symplectic, and Proposition \ref{posints}
guarantees that they will intersect each other positively, so
there will exist an almost complex structure making $C_s$
holomorphic. However, if $s_1$ and $s_2$ are two different
sections in the moduli space enumerated by
$\widetilde{\mathcal{DS}}_f(\alpha;\alpha_1,\ldots,\alpha_n)$, it
is unclear whether there will exist a single almost complex
structure on $X$ making both $C_{s_1}$ and $C_{s_2}$ holomorphic.

The almost complex structures on $X_r(f)$ used in the definition
of $\widetilde{\mathcal{DS}}$ are, quite crucially, required to
preserve the tangent space to the \emph{diagonal stratum}
consisting of divisors with one or more points repeated.  One
might hope to define analogous invariants which agree with
$Gr(\alpha;\alpha_1,\ldots,\alpha_n)$ using arbitrary almost
complex structures on $X_r(f)$.  If one could do this, though, the
arguments reviewed in Section \ref{review} would rather quickly
enable one to conclude that
$Gr(\alpha;\alpha_1,\ldots,\alpha_n)=0$ whenever $\alpha$ has
larger pairing with the symplectic form than does the canonical
class and $\alpha_i\cdot\alpha_j=0$ for $i\neq j$.  However, this
is not the case: the manifold considered in \cite{MT} admits a
symplectic form such that, for certain primitive, orthogonal,
square-zero classes $\alpha$, $\beta$, $\gamma$, and $\delta$ each
with positive symplectic area, the canonical class is
$2(\alpha+\beta+\gamma)$ but the invariant
$Gr(2(\alpha+\beta+\gamma)+\delta;\alpha,\beta,\gamma,\alpha+\beta+\gamma+\delta)$
is nonzero.

While the Gromov--Taubes invariant restricts attention to curves
whose components are all covers of embedded curves which do not
intersect each other, it is natural to hope for information about
curves Poincar\'e dual to $\alpha$ having some number $n$ of
transverse self-intersections. One might like to define an
analogue $Gr_n(\alpha)$ of the Gromov--Taubes invariant counting
such curves, but as we review in Section \ref{secfam}, owing to
issues relating to multiple covers it is somewhat unclear what the
definition of such an invariant should be in general.  If one
imposes some rather stringent conditions on $\alpha$ ($\alpha$
should be ``$n$-semisimple'' in the sense of Definition
\ref{sem}), there is however a natural such choice.

Note that for arbitrary $\alpha$ and $n$, following \cite{RT} one
may define an invariant $RT_{n}(\alpha)$ which might naively be
viewed as a count of \emph{connected} pseudoholomorphic curves
Poincar\'e dual to $\alpha$ with $n$ self-intersections by
enumerating solutions $u\co \Sigma_g\to X$ of the equation
$(\dbar_ju)=\nu (x,u(x))$ for generic $j$ and ``inhomogeneous term''
$\nu$, where the genus $g$ of the source curve is given in
accordance with the adjunction formula by
$2g-2=\alpha^2+\kappa_X\cdot\alpha-2n$.  (Note that the nontrivial
dependence of $\nu$ on $x$ prevents multiple cover problems from
arising.)  In the case $n=0$, the main theorem of \cite{IP2}
provides a universal formula equating $Gr(\alpha)$ with a certain
combination of the Ruan--Tian invariants $RT_0$.  The proof of that
theorem goes through easily to show that in the case when $\alpha$
is $n$-semisimple, there exists a similar formula equating
$Gr_n(\alpha)$ with a combination of Ruan--Tian invariants.  We
mention also that, again as an artifact of the multiple cover
problem, the Ruan--Tian invariants are obliged to take values in
$\mathbb{Q}$ rather than $\mathbb{Z}$.  $Gr(\alpha)$, on the other
hand, is an integer-valued invariant.

By combining the approaches of \cite{DS} and \cite{Liu}, in the
presence of a Lefschetz fibration $f\co X\to S^2$ we construct in
Section \ref{secfam} an integer-valued invariant
$\mathcal{FDS}_{f}^{n}(\alpha-2\sum e_i)$ which we conjecture to
be an appropriate candidate for a ``nodal version'' $Gr_n(\alpha)$
of the Gromov invariant for general classes $\alpha$ (after
perhaps dividing by $n!$ to account for a symmetry in the
construction). Pleasingly, the technical difficulties that often
arise in defining invariants like $Gr_n (\alpha)$ do not affect
$\mathcal{FDS}$: since $\mathcal{FDS}$ counts sections of a
(singular) fibration, which of course necessarily represent a
primitive homology class in the total space, we need not worry
about multiple covers; further, the fact that any bubbles that
form in the limit of a sequence of holomorphic sections must be
contained in the fibers of the fibration turns out (via an easy
elaboration of a dimension computation from \cite{DS}) to
generically rule out bubbling as well.  In principle, though,
$\mathcal{FDS}^{n}_{f}$ might depend on the choice of Lefschetz
fibration $f$.

Note that if $\pi\co X'\to X$ is a blowup with exceptional divisor
Poincar\'e dual to $\epsilon$, whenever $Gr_n(\beta)$ is defined
we will have $(Gr_n)_{X'}(\beta +\epsilon)=(Gr_n)_X(\beta)$ (here
and elsewhere we use the same notation for $\beta\in
H^2(X;\mathbb{Z})$ and $\pi^{*}\beta\in H^2(X';\mathbb{Z})$), as
the curves contributing to $(Gr_n)_X(\beta)$ generically miss the
point being blown up, and so the unions of their proper transforms
with the exceptional divisor will be precisely the curves
contributing to $(Gr_n)_{X'}(\beta +\epsilon)$.  With this said,
we formulate the:
\begin{conj} \label{conj}  Let $(X,\omega)$ be a symplectic 4-manifold and $\alpha\in H^2(X;\mathbb{Z})$, and
$f\co X'\to S^2$ a Lefschetz fibration obtained from a
sufficiently high-degree Lefschetz pencil on $X$, with the
exceptional divisors of the blowup $X'\to X$ Poincar\'e dual to
the classes $\epsilon_1,\ldots,\epsilon_N$.  Then the family
Donaldson--Smith invariants \[
\mathcal{FDS}^{n}_{f}\left(\alpha+\sum_{i=1}^{N}\epsilon_i-2\sum_{k=1}^{n}e_k\right)
\] are independent of the choice of $f$, and have a general
expression in terms of the Ruan--Tian invariants of $X$.
\end{conj}

Note that this conjectural general expression would then produce
an integer  by taking appropriate combinations of the (\emph{a
priori} only rational) Ruan--Tian invariants, similarly to the
formula of \cite{IP2}

In light of the behavior of $Gr_n$ under blowups, Theorem
\ref{famsame} amounts to the statement that:

\begin{thm} If $\alpha$ is strongly $n$-semisimple, then
Conjecture \ref{conj} holds for $\alpha$; more specifically, we
have
\[
\mathcal{FDS}^{n}_{f}\left(\alpha+\sum_{i=1}^{N}\epsilon_i-2\sum_{k=1}^{n}e_k\right)=n!Gr_n(\alpha)\] if $f$ has sufficiently high degree.
\end{thm}

We also prove that $\mathcal{FDS}$ vanishes under certain
circumstances.  This result depends heavily on the constructions
used by Smith in \cite{Smith} to prove his duality theorem, and so
we review these constructions in Section \ref{review}.  Section
\ref{vanproof} is then devoted to a proof of the following
theorem.

\begin{thm} \label{famvan} If $b^+(X)>b_1(X)+4n+1$, then for all $\alpha \in
H^2(X;\mathbb{Z})$ such that $r=\langle \alpha,[\Phi]\rangle$
satisfies $r>\max\{g(\Phi)+3n+d(\alpha),(4g(\Phi)-11)/3\}$, either
$\mathcal{FDS}^{n}_{f}(\alpha-2\sum e_i)=0$ or there exists an
almost complex structure $j$ on $X$ compatible with the fibration
$f\co X\to S^2$ which simultaneously admits holomorphic curves $C$
and $D$ Poincar\'e dual to the classes $\alpha$ and
$\kappa_X-\alpha$.  In particular,
$\mathcal{FDS}^{n}_{f}(\alpha-2\sum e_i)=0$ if $\alpha$ has larger
pairing with the symplectic form than does $\kappa_X$.
\end{thm}

Note that in the Lefschetz fibrations obtained from degree-$k$
Lefschetz pencils on some fixed symplectic manifold $(X,\omega)$,
the number $N$ of exceptional sections is $k^2[\omega]^2$ while the number $2g(\Phi)-2$
is asymptotic to $k^2[\omega]^2$, so the invariants
\[ \mathcal{FDS}^{n}_{f}\left(\alpha+\sum_{i=1}^{N}\epsilon_i-2\sum_k
e_k\right)\] considered in Conjecture \ref{conj} all eventually
satisfy the restriction on $r$ in Theorem \ref{famvan}.

The almost complex structure in the second alternative in Theorem
\ref{famvan} cannot be taken to be regular (in the sense that the
moduli spaces $\mathcal{M}^{j}_{X}(\beta)$ of $j$-holomorphic
curves Poincar\'e dual to $\beta$ are of the expected dimension);
the most we can say appears to be that it can be taken to be a
member of a regular $4n$-real-dimensional family of almost complex
structures, $i.e.$, a family of almost complex structures
$\{j_b\}$ parametrized by elements $b$ of an open set in
$\mathbb{R}^{4n}$ such that the spaces $\{(b,C)| C\in
\mathcal{M}^{j_b}_{X}(\beta)\}$ are of the expected real dimension
$2d(\beta)+4n$ near each $(b,C)$ such that $C$ has no
multiply-covered components. Also, if $X$ is in fact K\"ahler and
admits a compatible \emph{integrable} complex structure $j_0$ with
respect to which the fibration $f$ is holomorphic, then we can
take the $j$ in Theorem \ref{famvan} equal to $j_0$.

In fact, if we could take $j$ to be regular, then we could rule
out the second alternative in Theorem \ref{famvan} entirely (when
$n>0$) using the following argument: the invariant vanishes
trivially when $d(\alpha)< n$, so we can assume \\
$d(\alpha)=-\frac{1}{2}\alpha\cdot (\kappa-\alpha)>0$.  But then
our curves Poincar\'e dual to $\alpha$ and $\kappa-\alpha$ have
negative intersection number, which is only possible if they share
one or more components of negative square.  For generic $j$, a
virtual dimension computation shows that the only irreducible $j$-holomorphic
curves of negative square are $(-1)$-spheres.  Moreover whatever
$(-1)$-spheres appear in $X$ must be disjoint, since if they were
not, blowing one of two intersecting $(-1)$-spheres down would
cause the image of the other to be a symplectic sphere of
nonnegative self-intersection, which (by a result of \cite{McD})
would force $X$ to have $b^+=1$, which we assumed it did not.
Ignoring all the $(-1)$-spheres in $C$ and $D$ and taking the
union of what is left over would then give a $j$-holomorphic curve
Poincar\'e dual to a class $\kappa_X-\sum a_i e_i$ where the $e_i$
are classes of $(-1)$-spheres with $e_i\cdot e_k=0$ for $i\neq k$
and where at least one $a_i\geq 2$.  But one easily finds
$d(\kappa_X-\sum a_i e_i)<0$, so this too is impossible for
generic $j$.  For nongeneric $j$, this argument breaks down
because of the possibility that $C$ and $D$ might share components
of negative square and negative expected dimension, and there is a
wider diversity of possible homology classes of such curves.

The final section of the paper contains proofs of two technical
results that are used in the proofs of the main theorems. First,
we show that the operation of blowing up a point can be performed
in the almost complex category, a fact which does not seem to
appear in the literature and whose proof is perhaps more subtle
than one might anticipate.  The paper then closes with a proof of
the following result, which is necessary for the compactness
argument that we use to justify the definition of our invariant
$\mathcal{FDS}$: \begin{thm} Let $F\co \mathcal{H}_r\to D^2$
denote the $r$-fold relative Hilbert scheme of the map
$(z,w)\mapsto zw$, $\phi_0$ the partial resolution map
$F^{-1}(0)\to Sym^r\{zw=0\}$, and $\Delta\subset \mathcal{H}_r$
the diagonal stratum. At any point $p\in \Delta\cap F^{-1}(0)$
with $\phi_0(p)=\{(0,0),\ldots,(0,0)\}$, where $T_p\Delta$ is the
tangent cone to $\Delta$ at $p$, we have $T_p\Delta\subset T_p
F^{-1}(0)$.
\end{thm}

We end the introduction with some remarks on the possible relation
of $\mathcal{FDS}$ to (family) Seiberg--Witten theory. In
\cite{Sa} it was shown that where $X$ is the product of
$\mathbb{R}$ and a fibered three-manifold, so that $X$ fibers over
a cylinder, if one examines the Seiberg--Witten equations on $X$
using a family of metrics for which the size of the fibers shrinks
to zero, then one obtains in the adiabatic limit the equations for
a holomorphic family of solutions to the symplectic vortex
equations on the fibers.  In turn, there is a natural isomorphism
between the space of solutions to the vortex equations on a
Riemann surface and the symmetric product of the surface. In other
words, in this simple context the adiabatic limit of the
Seiberg--Witten equations is the equation for a holomorphic family
of elements of the symmetric products of the fibers of the
fibration $X\to \mathbb{R}\times S^1$.  As was noted in \cite{DS},
since for a Lefschetz fibration $f\co X\to S^2$ $\mathcal{DS}_f$
precisely counts pseudoholomorphic families of elements of the
symmetric products of the fibers of $f$, one might take
inspiration from Salamon's example and hope to obtain the
equivalence between $\mathcal{DS}_f$ and the Seiberg--Witten
invariant by considering the Seiberg--Witten equations on $X$ for
a family of metrics with respect to which the size of the fibers
shrinks to zero.

Now our invariant $\mathcal{FDS}^{n}_{f}$ is constructed by
counting pseudoholomorphic families of elements of the symmetric
products of the fibers of a family of Lefschetz fibrations $f^b$
obtained by restricting a map $f_n\co X_{n+1}\to S^2\times X_n$ to
the preimage $X^b$ of $S^{2}\times\{b\}$ as $b$ ranges over the
complement $X'_n$ of a set of codimension 4 in $X_n$. In the above
vein, one might hope to relate the family Seiberg--Witten
invariants $FSW$ for the family of $4$-manifolds $X_{n+1}\to X_n$
(which enumerate Seiberg--Witten monopoles in the various $X^b$ as
$b$ ranges over $X_n$; see, \emph{e.g.}, \cite{LL}) to
$\mathcal{FDS}^{n}_{f}$ via an adiabatic limit argument.  This
would in particular yield a proof of the independence of
$\mathcal{FDS}^{n}_{f}$ from $f$ in Conjecture \ref{conj}, and
indeed may well be the most promising way to establish this
independence in the absence of a suitable invariant $Gr_n$ (or of
a ``family Gromov--Taubes invariant'' $FGr$) with which
$\mathcal{FDS}^{n}_{f}$ might be directly equated.

As was shown in \cite{Liu2}, when $X$ is an algebraic surface and
$b^{+}(X)=1$ the family Seiberg--Witten invariants agree with
certain curve counts in algebraic geometry.  For larger values of
$b^+$, though, the family Seiberg--Witten invariants that are
hoped to correspond to nodal curve counts  are  expected to vanish
due to the fact that symplectic manifolds have Seiberg--Witten
simple type; note that Theorem \ref{famvan} suggests that
$\mathcal{FDS}^{n}_{f}$ also tends to vanish for large $b^+$.  By
contrast, there are plenty of  nontrivial nodal curve counts in
algebraic surfaces with $b^+>1$ (see \cite{Liu} for a review of
some of these); these counts correspond to Liu's ``algebraic
Seiberg--Witten invariants'' $\mathcal{ASW}$ and differ from $FSW$
when $b^+>1$.

\subsection*{Acknowledgements}  Section \ref{red} of this paper
appeared in my thesis \cite{thesis}; I would like to thank my
advisor Gang Tian for suggesting that I attempt to study nodal
curves using the Donaldson--Smith approach.  Thanks also to Cliff
Taubes for helping me identify an error in an earlier version of
this paper, to Dusa McDuff for making me aware of the need for
Section \ref{app2}, and to Ivan Smith for helpful remarks. This
work was partially supported by the National Science Foundation.

\section{Refining the standard surface
count} \label{red} Throughout this section, $X_r(f)$ will denote
the relative Hilbert scheme constructed from some high-degree but
fixed Lefschetz fibration $f\co X\to S^2$ obtained by Donaldson's
construction applied to the fixed symplectic 4-manifold
$(X,\omega)$.  The fiber of $f$ over $t\in S^2$ will occasionally
be denoted by $\Sigma_t$, and the homology class of the fiber by
$[\Phi]$.

As has been mentioned earlier, $\mathcal{DS}_f (\alpha)$ is a
count of holomorphic sections of the relative Hilbert scheme $X_r
(f)$ in a certain homotopy class $c_{\alpha}$ characterized by the
property that if $s$ is a section in the class $c_{\alpha}$ then
the closed set $C_s\subset X$ ``swept out'' by $s$ (that is, the
union over all $t$ of the divisors $s(t)\in\Sigma_t$) is
Poincar\'e dual to $\alpha$ (note that points of $C_s$ in this
interpretation may have multiplicity greater than 1).  That
$c_{\alpha}$ is the unique homotopy class with this property is
seen in Lemma 4.1 of \cite{Smith}; in particular, for instance, we
note that sections which descend to \emph{connected} standard
surfaces Poincar\'e dual to $\alpha$ are not distinguished at the
level of homotopy from those which descend to disjoint unions of
several standard surfaces which combine to represent $PD(\alpha)$.

Of course, in studying standard surfaces it is natural to wish to
know their connected component decompositions, so we will
presently attempt to shed light on this.  Suppose that we have a
decomposition \[ \alpha=\alpha_1 +\cdots +\alpha_n \] with
\[ \langle \alpha,[\Phi]\rangle=r, \quad \langle
\alpha_i,[\Phi]\rangle =r_i. \]  Over each $t\in S^2$ we have an
obvious ``divisor addition map'' \begin{align}
+\co\prod_{i=1}^{n}S^{r_i}\Sigma_t&\to S^r\Sigma_t \nonumber \\
(D_1,\ldots,D_n)&\mapsto D_1+\cdots+D_n; \nonumber \end{align}
allowing $t$ to vary we obtain from this a map on sections:
\begin{align} +\co\prod_{i=1}^{n}\Gamma(X_{r_i}(f))&\to \Gamma(X_r (f))\nonumber \\
(s_1,\ldots,s_n)&\mapsto \sum_{i=1}^n s_i. \nonumber \end{align}
As should be clear, one has \[ +(c_{\alpha_1}\times\cdots\times
c_{\alpha_n})\subset c_\alpha\] if $\alpha =\sum \alpha_i$, since
$C_{\sum \alpha_i}$ is the union of the standard surfaces
$C_{s_i}$ and hence is Poincar\'e dual to $\alpha$ if each
$C_{s_i}$ is Poincar\'e dual to $\alpha_i$.  Further, we readily
observe:
\begin{lemma} \label{closed} The image
$+(c_{\alpha_1}\times\cdots\times c_{\alpha_n})\subset c_{\alpha}$
is closed with respect to the $C^0$ norm.
\end{lemma}
\begin{proof}
Suppose we have a sequence
$(s^{m}_{1},\ldots,s^{m}_{n})_{m=1}^{\infty}$ in
$c_{\alpha_1}\times\cdots\times c_{\alpha_n}$ such that $\sum
s_{i}^{m}\to s\in c_{\alpha}$.  Now each $S^{r_i}\Sigma_t$ is
compact, so at each $t$, each of the sequences $s_{i}^{m}(t)$ must
have subsequences converging to some $s_{i}^{0}(t)$. But then
necessarily each $\sum s_{i}^{0}(t)=s(t)$, and then we can see by,
for any $l$, fixing the subsequence used for all $i\neq l$ and
varying that used for $i=l$ that in fact every subsequence of
$s_{l}^{m}(t)$ must converge to $s_{l}^{0}(t)$.  Letting $t$ vary
then gives sections $s_{i}^{0}$ such that every $s_{i}^{m}\to
s_{i}^{0}$ and $\sum s_{i}^{0}=s$; the continuity of $s$ is
readily seen to imply that of the $s_{i}^{0}$.
\end{proof}

At this point it is useful to record an elementary fact about the
linearization of the divisor addition map.

\begin{prop} \label{plus} Let $\Sigma$ be a Riemann surface and $r=\sum r_i$.  The linearization $+_*$ of the addition map
\[ +\co\prod_{i=1}^{n}S^{r_i}\Sigma\to S^r\Sigma\] at $(D_1,\ldots ,D_n)$ is an isomorphism if and only if
$D_i\cap D_j=\varnothing$ for $i\neq j$.  If two or more of the
$D_i$ have a point in common, then the image of $+_*$ at
$(D_1,\ldots,D_n)$ is contained in $T_{\sum D_i}\Delta$, where
$\Delta\subset S^r\Sigma$ is the diagonal stratum consisting of
divisors with a repeated point.
\end{prop}
\begin{proof}
By factoring $+$ as a composition \[ S^{r_1}\Sigma\times
S^{r_2}\Sigma\times\cdots\times S^{r_n}\Sigma\to
S^{r_1+r_2}\Sigma\times\cdots\times S^{r_n}\Sigma\to\cdots\to
S^{r}\Sigma \] in the obvious way we reduce to the case $n=2$.
Now in general for a divisor $D=\sum a_i p_i \in S^d\Sigma$ where
the $p_i$ are distinct, a chart for $S^d \Sigma$ is given by
$\prod S^{a_i}U_i$, where the $U_i$ are holomorphic coordinate
charts around $p_i$ and the $S^{a_i}U_i$ use as coordinates the
elementary symmetric polynomials $\sigma_1,\ldots,\sigma_{a_i}$ in
the coordinates of $U_{i}^{a_i}$.  As such, if $D_1$ and $D_2$ are
disjoint, a chart around $D_1 +D_2 \in S^{r_1 +r_2}\Sigma$ is
simply the Cartesian product of charts around $D_1\in
S^{r_1}\Sigma$ and $D_2\in S^{r_2}\Sigma$, and the map $+$ takes
the latter diffeomorphically (indeed, biholomorphically) onto the
former, so that $(+_*)_{(D_1,D_2)}$ is an isomorphism.

On the other hand, note that \[ +\co S^a\mathbb{C}\times
S^b\mathbb{C}\to S^{a+b}\mathbb{C} \] is given in terms of the
local elementary symmetric polynomial coordinates around the
origin by \[
(\sigma_1,\ldots,\sigma_a,\tau_1,\ldots,\tau_b)\mapsto
(\sigma_1+\tau_1,\sigma_2+\sigma_1\tau_1+\tau_2,\ldots,\sigma_a\tau_b),\]
and so has linearization \[
(+_*)_{(\sigma_1,\ldots,\tau_b)}(\eta_1,\ldots,\eta_a,\zeta_1,\ldots,\zeta_b)=
(\eta_1+\zeta_1,\eta_2+\sigma_1\zeta_1+\tau_1\eta_1+\zeta_2,\ldots,\sigma_a\zeta_b+\tau_b\eta_a).\]
We thus see that $Im(+_*)_{(0,\ldots,0)}$ only has dimension
$\max\{a,b\}$ and is contained in the image of the linearization
of the smooth model
\begin{align} \mathbb{C}\times S^{a+b-2}\mathbb{C}&\to S^{a+b}\mathbb{C} \nonumber \\ (z,D)&\mapsto 2z+D \nonumber \end{align}
for the diagonal stratum at $(0,0+\cdots +0)$. Suppose now that
$D_1$ and $D_2$ contain a common point $p$; write $D_i =a_i p+D_i
'$ where $D_i\in S^{r_i-a_i}\Sigma$ are divisors which do not
contain $p$.  Then from the commutative diagram \[
\begin{CD}
{S^{a_1}\Sigma\times S^{r_1-a_1}\Sigma\times S^{a_2}\Sigma\times S^{r_2-a_2}\Sigma}@>>>{S^{r_1}\Sigma\times S^{r_2}\Sigma}\\
@VVV @ VV{+}V\\
{S^{a_1+a_2}\Sigma\times S^{r_1+r_2-a_1-a_2}\Sigma}@>>>{S^{r_1+r_2}\Sigma}\\
\end{CD}\] and the fact that the linearization of the top arrow at $(a_1 p,D_{1}',a_2 p,D_{2}')$ is an isomorphism (by what we
showed earlier, since the $D_{i}'$ do not contain $p$), while the
linearization of the composition of the left and bottom arrows at
$(a_1 p,D_{1}',a_2 p,D_{2}')$ has image contained in
$T_{D_1+D_2}\Delta$, it follows that $(+_*)_{(D_1,D_2)}$ has image
contained in $T_{D_1+D_2}\Delta$ as well, which suffices to prove
the proposition.
\end{proof}

\begin{cor} \label{tgt} If $s_i \in \Gamma (X_{r_i}(f))$ are differentiable sections such that
$C_{s_i}\cap C_{s_j} \neq\varnothing$ for some $i\neq j$, then
$s=\sum s_i \in \Gamma(X_r(f))$ is tangent to the diagonal stratum
of $X_r (f)$.  \end{cor}
\begin{proof}
Indeed, if $C_{s_i}\cap C_{s_j} \neq\varnothing$, then there is
$x\in S^2$ such that the divisors $s_i (x)$ and $s_j (x)$ contain
a point in common, and so for $v\in T_x S^2$ we have \[ s_*
v=(+\circ(s_i,s_j))_* v=+_*(s_{1*}v,s_{2*}v)\in T_{s(t)}\Delta
\] by Proposition \ref{plus}. \end{proof}

Note that it is straightforward to find cases in which the $s_i$
are only continuous with some $C_{s_i}\cap C_{s_j}$ nonempty and
the sum $s=\sum s_i$ is smooth but not tangent to the diagonal.
For example, let $r=2$, and in local coordinates let $s_1$ be a
square root of the function $z\mapsto Re(z)$ and $s_2=-s_1$.  Then
in the standard coordinates on the symmetric product we have
$s(z)=(0,-Re(z))$, so that $T(Im s)$ shares only one dimension
with $T\Delta$ at $z=0$.  If $s$ is \emph{transverse} to $\Delta$,
one can easily check that a similar situation cannot arise.

We now bring pseudoholomorphicity in the picture. Throughout this
treatment, all almost complex structures on $X_r (f)$ will be
assumed to agree with the standard structures on the symmetric
product fibers, to make the map $F\co X_r (f)\to S^2$
pseudoholomorphic, and, on some (\emph{not} fixed) neighborhood of
the critical fibers of $F$, to agree with the holomorphic model
for the relative Hilbert scheme over a disc around a critical
value for $f$ provided in Section 3 of \cite{Smith}.  Let
$\mathcal{J}$ denote the space of these almost complex structures.
It follows by standard arguments (see Proposition 3.4.1 of
\cite{MS} for the general scheme of these arguments and Section 4
of \cite{DS} for their application in the present context) that
for generic $J\in\mathcal{J}$ the space
$\mathcal{M}^{J}(c_{\alpha})$ is a smooth manifold of (real)
dimension $2d(\alpha)=\alpha^2-\kappa_{X}\cdot\alpha$ (the
dimension computation comprises Lemma 4.3 of \cite{Smith}); this
manifold is compact, for bubbling is precluded by the arguments of
Section 4 of \cite{Smith} assuming we have taken a sufficiently
high-degree Lefschetz fibration.

Inside $\mathcal{M}^{J}(c_{\alpha})$ we have the set
$\mathcal{M}^{J}(c_{\alpha_1}\times\cdots\times c_{\alpha_n})$
consisting of holomorphic sections which lie in the image
$+(c_{\alpha_1}\times\cdots\times c_{\alpha_n})$.  By Lemma
\ref{closed} and the compactness of $\mathcal{M}^{J}(c_{\alpha})$,
$\mathcal{M}^{J}(c_{\alpha_1}\times\cdots\times c_{\alpha_n})$ is
evidently compact; however, the question of its dimension or even
whether it is a manifold appears to be a more subtle issue in
general.

Let us pause to consider what we would like the dimension of
$\mathcal{M}^{J}(c_{\alpha_1}\times\cdots\times c_{\alpha_n})$ to
be.  The objects in
$\mathcal{M}^{J}(c_{\alpha_1}\times\cdots\times c_{\alpha_n})$ are
expected to correspond in some way to unions of holomorphic curves
Poincar\'e dual to $\alpha_i$.  Accordingly, assume we have chosen
the $\alpha_i$ so that
$d(\alpha_i)=\frac{1}{2}(\alpha_{i}^{2}-\kappa_{X}\cdot\alpha_i)\geq
0$ (for otherwise we would expect
$\mathcal{M}^{J}(c_{\alpha_1}\times\cdots\times c_{\alpha_n})$ to
be empty).  Holomorphic curves in these classes will intersect
positively as long as they do not share any components of negative
square; for a generic almost complex structure the only such
components that can arise are $(-1)$-spheres, so if we choose the
$\alpha_i$ to not share any $(-1)$-sphere components (\emph{i.e.},
if the $\alpha_i$ are chosen so that there is no class $E$
represented by a symplectic $(-1)$-sphere such that $\langle
\alpha_i , E\rangle<0$ for more than one $\alpha_i$), then it
would also be sensible to assume that $\alpha_i\cdot\alpha_j\geq
0$ for $i\neq j$.

The above naive interpretation of
$\mathcal{M}^{J}(c_{\alpha_1}\times\cdots\times c_{\alpha_n})$
would suggest that its dimension ought to be $\sum d(\alpha_i)$.
Note that \[ d(\alpha)=d(\sum\alpha_i)=\sum
d(\alpha_i)+\sum_{i>j}\alpha_i\cdot\alpha_j, \] so under the
assumptions on the $\alpha_i$ from the last paragraph we have that
the expected dimension of
$\mathcal{M}^{J}(c_{\alpha_1}\times\cdots\times c_{\alpha_n})$ is
at most the actual dimension of $\mathcal{M}^J(c_{\alpha})$ (as we
would hope, given that the former is a subset of the latter), with
equality if and only if $\alpha_i\cdot\alpha_j=0$ whenever $i\neq
j$.

As usual, we will find it convenient to cut down the dimensions of
our moduli spaces by imposing incidence conditions,
 so we shall fix a set $\Omega$ of points $z\in X$ and consider the space
$\mathcal{M}^{J,\Omega}(c_{\alpha_1}\times\cdots\times
c_{\alpha_n})$ of elements
$s\in\mathcal{M}^{J}(c_{\alpha_1}\times\cdots\times c_{\alpha_n})$
such that $C_s$ passes through each of the points $z$ (or, working
more explicitly in $X_r(f)$, such that $s$ meets each divisor
$z+S^{r-1}\Sigma_t$, $\Sigma_t$ being the fiber which contains
$z$).  $\mathcal{M}^{J,\Omega}(c_{\alpha})$ is defined similarly,
and standard arguments show that for generic choices of $\Omega$
$\mathcal{M}^{J,\Omega}(c_{\alpha})$ will be a compact manifold of
dimension \[ 2(d(\alpha)-\#\Omega).\]

We wish to count $J$-holomorphic sections $s$ of $X_r(f)$ such
that the reducible components of $C_s$ are Poincar\'e dual to the
$\alpha_i$. If we impose $\sum d(\alpha_i)$ incidence conditions,
then according to the above discussion
$\mathcal{M}^{J,\Omega}(c_{\alpha})$ will be a smooth manifold of
dimension $2\sum_{i>j}\alpha_i\cdot\alpha_j$. A section $\sum
s_i\in +(c_{\alpha_1}\times\cdots c_{\alpha_n})$ whose summands
are all differentiable would then, by Corollary \ref{tgt}, have
one tangency to the diagonal $\Delta$ for each of the
intersections between the $C_{s_i}$, of which the total expected
number is $\sum_{i>j}\alpha_i\cdot\alpha_j$.  This suggests that
the sections we wish to count should be found among those elements
of $\mathcal{M}^{J,\Omega}(c_{\alpha})$ which have
$\sum_{i>j}\alpha_i\cdot\alpha_j$ tangencies to $\Delta$, where
$\Omega$ is a set of $\sum d(\alpha_i)$ points.

To count pseudoholomorphic curves tangent to a symplectic
subvariety it is necessary to restrict to almost complex
structures which preserve the tangent space to the subvariety (see
\cite{IP} for the general theory when the subvariety is a
submanifold).  Accordingly, we shall restrict attention to the
class of almost complex structures $J$ on $X_r (f)$ which are
\emph{compatible with the strata} in the sense to be explained
presently (for more details, see Section 6 of \cite{DS}, in which
the notion was introduced).

Within $\Delta$, there are various strata $\chi_{\pi}$ indexed by
partitions $\pi: r=\sum a_i n_i$ with at least one $a_i >1$; these
strata are the images of the maps \begin{align}
p_{\chi}\co X_{n_1}(f)\times_{S^2}\cdots\times_{S^2}X_{n_k}(f)&\to X_r (f) \nonumber\\
(D_1,\ldots ,D_k)&\mapsto \sum a_i D_i; \nonumber \end{align}
 in particular, $\Delta = \chi_{r=2\cdot 1+1\cdot (r-2)}$.  An almost complex structure $J$ on $X_r (f)$ is said to be compatible with the strata if the maps $p_{\chi}$ are
$(J',J)$-holomorphic for suitable almost complex structures $J'$
on their domains.

Denoting by $Y_{\chi}$ the domain of $p_{\chi}$, Lemma 7.4 of
\cite{DS} and the discussion preceding it show:

\begin{lemma}[\cite{DS}]\label{DS}  For almost complex structures $J$ on $X_r (f)$ which are compatible with the strata, each $J$-holomorphic section $s$ of $X_r (f)$ lies in some
unique minimal stratum $\chi$ and meets all strata contained in
$\chi$ in isolated points.  In this case, there is a
$J'$-holomorphic section $s'$ of $Y_{\chi}$ such that
$s=p_{\chi}\circ s'$. Furthermore, for generic $J$ among those
compatible with the strata, the actual dimension of the space of
all such sections $s$ is equal to the expected dimension of the
space of $J'$-holomorphic sections $s'$ lying over $s$.
\end{lemma}

We note the following analogue for standard surfaces of the
positivity of intersections of pseudoholomorphic curves.

\begin{prop} \label{posints}
Let $s=m_1 s_1+\cdots +m_k s_k$ be a $J$-holomorphic section of
$X_r (f)$, where the $s_i\in c_{\alpha_i}\subset \Gamma
(X_{r_i}(f))$ are each not contained in the diagonal stratum of
$X_{r_i}(f)$, and where the almost complex structure $J$ on $X_r
(f)$ is compatible with the strata.  Assume that the $s_i$ are all
differentiable.  Then all isolated intersection points of
$C_{s_i}$ and $C_{s_j}$ contribute positively to the intersection
number $\alpha_i\cdot \alpha_j$.
\end{prop}
\begin{proof}
We shall prove the lemma for the case $k=2$, the general case
being only notationally more complicated.  The analysis is
somewhat easier if the points of $C_{s_1}\cap C_{s_2}\subset X$ at
issue only lie over $t\in S^2$ for which $s_1 (t)$ and $s_2 (t)$
both miss the diagonal of $X_{r_1}(f)$ and $X_{r_2}(f)$,
respectively, so we first argue that we can reduce to this case.
Let $\chi$ be the minimal stratum (possibly all of $X_r (f)$) in
which $s=m_1s_1+m_2s_2$ is contained, so that all intersections of
$s$ with lower strata are isolated.  Let $p\in X$ be an isolated
intersection point of $C_{s_1}$ and $C_{s_2}$ lying over $0\in
S^2$, and let $\delta
>0$ be small enough that there are  no other intersections of $s$ with any substrata of $\chi$ (and so in particular no other points of $C_{s_1}\cap C_{s_2}$) lying
over $D_{2\delta}(0)\subset S^2$. We may then perturb
$s=m_1s_1+m_2s_2$ to $\tilde{s}=m_1\tilde{s_1}+m_2\tilde{s_2}$,
still lying in $\chi$, such that\\ \begin{itemize} \item[(i)] Over
$D_{\delta}(0)$, $\tilde{s}$ is $J$-holomorphic and disjoint from
all substrata having real codimension larger than 2 in $\chi$, and
the divisors $\tilde{s_1}(0)$ and $\tilde{s_2}(0)$ both still
contain $p$; \item[(ii)] Over the complement of $D_{2\delta}(0)$,
$\tilde{s}$ agrees with $s$; and \item[(iii)] Over
$D_{2\delta}(0)\setminus D_{\delta}(0)$, $\tilde{s}$ need not be
$J$-holomorphic but is connected to $s$ by a family of sections
$s_t$ contained in $\chi$ which miss all substrata of $\chi$
\end{itemize} (it may be necessary to decrease $\delta$ to find
such $\tilde{s}$, but after doing so such $\tilde{s}$ will exist
by virtue of the abundance of $J$-holomorphic sections over the
small disc $D_{\delta}(0)$ which are close to
$s|_{D_{\delta}(0)}$).  The contribution of $p$ to the
intersection number $\alpha_1\cdot\alpha_2$ will then be equal to
the total contribution of all the intersections of
$C_{\tilde{s_1}}$ and $C_{\tilde{s_2}}$ lying over
$D_{\delta}(0)$, and the fact that $\tilde{s}$ misses all
substrata with codimension larger than 2 in $\chi$ is easily seen
to imply that these intersections (of which there is at least one,
at $p$) are all at points where $\tilde{s_1}$ and $\tilde{s_2}$
miss the diagonals in $X_{r_1}(f)$ and $X_{r_2}(f)$.

As such, it suffices to prove the lemma for intersection points at
which $s_1$ and $s_2$ both miss the diagonal.  In this case, in a
coordinate neighborhood $U$ around $p$, the $C_{s_i}$ can be
written as graphs $C_{s_i}\cap U=\{w=g_i (z)\}$, where $w$ is the
holomorphic coordinate on the fibers of $X$, $z$ is the pullback
of the holomorphic coordinate on $S^2$, and $g_i$ is a
differentiable complex-valued function which vanishes at $z=0$.
Suppose first that $m_1=m_2=1$.  Then near $s(0)$, we may use
coordinates $(z,\sigma_1,\sigma_2,y_3,\ldots,y_r)$ for $X_r (f)$
obtained from the splitting $T_0S^2\oplus T_{2p}S^2\Sigma_0\oplus
T_{s(t)-2p}S^{r-2}\Sigma_0$, and the first two vertical
coordinates of $s(z)=(s_1+s_2)(z)$ with respect to this splitting
are $(g_1(z)+g_2(z),g_1(z)g_2(z))$.  Now $s$ is $J$-holomorphic
and meets the $J$-holomorphic diagonal stratum $\Delta$ at
$(0,s(0))$, and at this point $\Delta$ is tangent to the
hyperplane $\sigma_2=0$, so it follows from Lemma 3.4 of \cite{IP}
that the Taylor expansion of $g_1 (z)g_2(z)$ has form $a_0 z^{d}
+O(d+1)$.  But then the Taylor expansions of $g_1 (z)$ and $g_2
(z)$ begin, respectively, $a_1 z^{d_1}+O(d_1+1)$ and $a_2
z^{d_2}+O(d_2+1)$, with $d_1+d_2=d$.  Then since $C_{s_i}\cap
U=\{w=g_i (z)\}$, it follows immediately that the $C_{s_i}$ have
intersection multiplicity $\max\{d_1,d_2\}>0$ at $p$.

There remains the case where one or both of the $m_i$ is larger
than 1.  In this case, where
$Y_{\chi}=X_{r_1}(f)\times_{S^2}X_{r_2}(f)$ is the smooth model
for $\chi$, because $J$ is compatible with the strata, $(s_1,s_2)$
is a $J'$-holomorphic section of $Y_{\chi}$ for an almost complex
structure $J'$ such that $p_{\chi}\co Y_{\chi}\to X_r (f)$ is
$(J',J)$-holomorphic.  Now where $\tilde{\Delta}=\{(D_1,D_2)\in
Y_{\chi}|D_1\cap D_2\neq \varnothing\}$, compatibility with the
strata implies that $\tilde{\Delta}$ will be $J'$-holomorphic. In
a neighborhood $V$ around $(s_1 (z),s_2(z))$, we have, in
appropriate coordinates, $\tilde{\Delta}\cap
V=\{(z,w,w,D_1,D_2)|w\in \Sigma_z\}$,while $(s_1 (z),s_2 (z))$ has
first three coordinates $(z,g_1 (z),g_2 (z))$.  From this it
follows by Lemma 3.4 of \cite{IP} that \[ g_1 (z)-g_2 (z)=a_0 z^d
+O(d+1) \] for some $d$, in which case $C_{s_1}$ and $C_{s_2}$
have intersection multiplicity $d>0$ at $p$.

\end{proof}

\begin{dfn} Let $\Omega$ be a set of $\sum d(\alpha_i)$ points and
let $J$ be an almost complex structure compatible with the strata.
$\mathcal{M}^{J,\Omega}_{0}(\alpha_1,\ldots,\alpha_n)$ shall
denote the set of $J$-holomorphic sections $s\in c_{\alpha}$ with
$\Omega\subset C_s$ such that there exist $C^1$ sections $s_i\in
c_{\alpha_i}$ with $s=\sum s_i$, while the $s_i$ themselves do not
admit nontrivial decompositions as sums of $C^1$ sections.
\end{dfn}

 We would like to assert that for generic $J$ and $\Omega$, the space
$\mathcal{M}^{J,\Omega}_{0}(\alpha_1,\ldots,\alpha_n)$ does not
include any sections contained within the strata.  This is not
true in full generality; rather we need the following assumption
in order to rule out the effects of multiple covers of square-zero
tori and $(-1)$-spheres in $X$.

\begin{ass} \label{ass} None of the $\alpha_i$ can be written as $\alpha_i =m\beta$ where
$m>1$ and either $\beta^2=\kappa_{X}\cdot\beta=0$ or
$\beta^2=\kappa_X\cdot\beta=-1.$\end{ass}

Under this assumption, we note that if $s=\sum s_i\in
\mathcal{M}^{J,\Omega}_{0}(\alpha_1,\ldots,\alpha_n)$ were
contained in $\Delta$, then since the $\alpha_i$ and hence the
$s_i$ are distinct we can write each $s_i$ as $s_i =
m_i\tilde{s_i}$ with at least one $m_i >1$.  The minimal stratum
of $s$ will then be $\chi_{\pi}$ where $\pi = \left\{r=\sum
m_i\left(\frac{r_i}{m_i}\right)\right\}$ and
$s'=(\tilde{s_1},\ldots,\tilde{s_n})$ will be a $J'$-holomorphic
section of $Y_{\chi}$ with $s=p_{\chi}\circ s'$, in the homotopy
class $[c_{\alpha_1/m_1}\times\cdots\times c_{\alpha_n/m_n}]$.

If any of the $d(\alpha_i/m_i)<0$, then Lemma \ref{DS} implies
that there will be no such sections $s'$ at all; otherwise (again
by Lemma \ref{DS}) the real dimension of the space of such
sections (taking into account the incidence conditions) will be
\begin{equation}\label{dim} 2\left(\sum d(\alpha_i/m_i)-\sum d(\alpha_i)\right). \end{equation}
  But an easy manipulation of the general formula for $d(\beta)$ and the adjunction formula (which applies here because the
standard surface corresponding to a section of $X_r (f)$ which
meets $\Delta$ positively will be symplectic; c.f. Lemma 2.8 of
\cite{DS}) shows that if $d(\beta)\geq 0$ and $m\geq 2$ then
$d(m\beta)> d(\beta)$  unless either
$\beta^2=\kappa_{X}\cdot\beta=0$ or
$\beta^2=\kappa_{X}\cdot\beta=-1$, and these are ruled out in this
context by (i) and (ii) above, respectively.  So Assumption
\ref{ass} implies that the dimension in Equation \ref{dim} is
negative, so no such $s'$ will exist for generic $J$.  This proves
part of the following:

\begin{prop} \label{cpct} Under Assumption \ref{ass}, for generic pairs $(J,\Omega)$ where $J$ is compatible with the strata and $\#\Omega =\sum d(\alpha_i)$,
$\mathcal{M}^{J,\Omega}_{0}(\alpha_1,\ldots,\alpha_n)$ is a finite
set consisting only of sections not contained in $\Delta$.
\end{prop}

\begin{proof}
That no member of
$\mathcal{M}^{J,\Omega}_{0}(\alpha_1,\ldots,\alpha_n)$ is
contained in $\Delta$ follows from the above discussion. As for
the dimension of our moduli space, note that any $s=\sum
s_i\in\mathcal{M}^{J,\Omega}_{0}(\alpha_1,\ldots,\alpha_n)$ has
one tangency (counted with multiplicity) to $\Delta$ for each of
the intersections of the $C_{s_i}$, of which there are $\sum
\alpha_i\cdot\alpha_j$ (counted with multiplicity; this
multiplicity will always be positive by Proposition
\ref{posints}). By the results of section 6 of \cite{IP}, the
space $\mathcal{M}^{J,\Omega}_{\delta,\Delta}(c_{\alpha})$ of
$J$-holomorphic sections in the class $c_{\alpha}$ having $\delta$
tangencies to $\Delta$ and whose descendant surfaces pass through
$\Omega$ will, for generic $(J,\Omega)$, be a manifold of
dimension \[ 2(d(\alpha)-\sum d(\alpha_i)-\delta)=2(\sum
\alpha_i\cdot \alpha_j-\delta), \] which is equal to zero in the
case $\delta=\sum \alpha_i\cdot\alpha_j$ of present relevance to
us.

Let us now show that
$\mathcal{M}^{J,\Omega}_{0}(\alpha_1,\ldots,\alpha_n)$ is compact.
Now since $+(c_{\alpha_1}\times\cdots \times c_{\alpha_n})$ is
$C^0$-closed in $c_{\alpha}$, by Gromov compactness any sequence
$s^{(m)}=\sum_{i=1}^{n} s^{(m)}_{i}$ in
$\mathcal{M}^{J,\Omega}_{0}(\alpha_1,\ldots,\alpha_n)$ has (after
passing to a subsequence) a $J$-holomorphic limit $s=\sum s_i$
where the $s_i\in c_{\alpha_i}$ are at least continuous.  We claim
that, at least for generic $(J,\Omega)$, we can guarantee the
$s_i$ to be $C^1$.  In light of Proposition \ref{plus}, the
differentiability of the $s_i$ is obvious at all points where $s$
misses the diagonal, since $s$ is smooth by elliptic regularity
and the divisor addition map induces an isomorphism on the tangent
spaces away from the diagonal.  Now each $s^{(m)}$ has $\sum
\alpha_i\cdot\alpha_j$ tangencies to the diagonal, corresponding
to points $t\in S^2$ at which some pair of the divisors
$s^{(m)}_{i}(t)$ share a point in common. The limit $s$ will then
likewise have $n$ tangencies to the diagonal; the dimension
formulas in \cite{IP} ensure that for generic $(J,\Omega)$ no two
of the tangencies will coalesce into a higher order tangency to
the smooth part of $\Delta$ in the limit, and all of the
intersecions of $Im\, s$ with the smooth part of the diagonal
other than these $n$ tangencies will be transverse. Furthermore,
one may easily show (using for instance an argument similar to the
one used in Lemma 2.1 of \cite{Usher} to preclude generic
0-dimensional moduli spaces of pseudoholomorphic curves in a
Lefschetz fibration from meeting the critical points) that since
the singular locus of $\Delta$ has codimension 4 in $X_r(f)$, if
$J$ has been chosen generically then $s$ will not meet
$\Delta^{sing}$, and so no $s(t)$ will contain more than one
repeated point (and that point cannot appear with multiplicity
larger than two).  In light of this, each tangency of $s$ to
$\Delta$ will occur at a point $s(t)$ where some pair $s_i(t)$ and
$s_j(t)$ have some point $p$ in common, and all other points
contained in any $s_k(t)$ are distinct from each other and from
$p$. Thanks to Proposition \ref{plus}, this effectively reduces us
to the case $r=2$, with $s=s_1+s_2$ a sum of continuous sections
with $s_1(0)=s_2(0)=0$ which is holomorphic with respect to an
almost complex structure which preserves the diagonal stratum
$\Delta$ in $D^2\times Sym^2D^2$, such that $s$ is tangent to
$\Delta$.  Then letting $\delta(t)=(s_1+s_2)^2(t)-4s_1(z)s_2(t)$
be the discriminant, that $s$ is tangent to the diagonal stratum
implies, using Lemma 3.4 of \cite{IP}, that $\delta(t)=at^2+O(3)$
for some constant $a$; in particular $\delta(t)$ has two $C^1$
square roots $\pm r(t)$. Since $s$ is smooth, so is its first
coordinate $t\mapsto s_1(t)+s_2(t)$; adding this smooth function
to the $C^1$ functions $\pm r(t)$ and dividing by two then
recovers the functions $s_1(t)$ and $s_2(t)$ and verifies that
they are $C^1$ at $t=0$.

We have thus shown that the $s_i$ are all $C^1$ at the points
where $s=\sum s_i$ is tangent to $\Delta$.  Where $s$ is
transverse to $\Delta$, one sees easily that the $s_i$ are
pairwise disjoint, with one $s_i$ transverse to the diagonal in
$X_{r_i}(f)$ and all others missing their diagonals, so the
differentiability of the $s_i$ is clear. This indeed verifies that
the limit $s=\sum s_i$ is a sum of $C^1$ sections $s_i$, since our
generic choice of $J$ is such that the only intersections of $Im
\, s$ with $\Delta$ only are either transverse or of second order.

Now each of the $C_{s^{(m)}_{i}}$ is connected, so $C_{s_i}$ is
connected as well. \emph{A priori}, it is possible that $s$ might
not lie in $\mathcal{M}^{J,\Omega}_{0}(\alpha_1,\ldots,\alpha_n)$
because some of the $s_i$ might decompose further, say as $s_i=m_1
u_{i_1}+\cdots+m_lu_{i_l}$ where $u_{i_j}\in c_{\beta_{i_j}}$ are
$C^1$. But since $C_{s_i}$ is connected, the $C_{u_{i_j}}$ cannot
all be disjoint, and by Corollary \ref{tgt} any intersection
between two of them would give rise to an additional tangency of
$s$ to $\Delta$, over and above the $n$ tangencies arising from
the intersections between distinct $C_{s_i}$.  Once again, this is
ruled out for generic $J$ by the dimension formulas of \cite{IP}.
%and so using Lemma \ref{posints} we have \[
%d(\alpha_i)=d(\sum m_j \beta_{i_j})=\sum d(m_j\beta_{i_j})+\sum
%m_j m_k \beta_{i_j}\cdot\beta_{i_k} > \sum d(m_j \beta_{i_j}).\]
%Meanwhile, Lemma \ref{DS} allows us to rule this decomposition out
%for generic $(J,\Omega)$ as long as \[ d(\alpha_i)>\sum
%d(\beta_{i_j}). \]
%   Now as noted earlier we have
%$d(m_j \beta_{i_j})>d(\beta_{i_j})$ unless $\beta_{i_j}$ is the
%class either of a $(-1)$ sphere or of square-zero torus; in the
%latter case we have $d(m_j\beta_{i_j})=d(\beta_{i_j})=0$, so it
%suffices to rule out decompositions of form
%$\alpha_i=(\alpha_i-mE)+mE$ where $E$ is the class of a $(-1)$
%sphere.  Such a decomposition will not occur if $\alpha_i=E$ (for
%then $C_{\alpha_i-mE}$ would not be symplectic), so by Assumption
%\ref{ass}(ii) we have $\alpha_i\cdot E\geq 0$ and so
%\[
%d(\alpha_i-mE)+d(E)=d(\alpha_i-mE)=d(\alpha)+d(-mE)+\alpha\cdot(-mE)\]\[\leq
%d(\alpha)+d(-mE)=d(\alpha)-(m^2+m)/2<d(\alpha).\]
This proves
that (for generic $J$) the summands $s_i$ in a sequence $s=\sum
s_i$ occurring as a limit point of
$\mathcal{M}^{J,\Omega}_{0}(\alpha_1,\ldots,\alpha_n)$ cannot
decompose further and hence themselves lie in
$\mathcal{M}^{J,\Omega}_{0}(\alpha_1,\ldots,\alpha_n)$, so that
$\mathcal{M}^{J,\Omega}_{0}(\alpha_1,\ldots,\alpha_n)$ is compact.

 Since we have already shown that
$\mathcal{M}^{J,\Omega}_{0}(\alpha_1,\ldots,\alpha_n)$ is
zero-dimensional, the proposition follows.

\end{proof}

\begin{prop} For generic $(J_0,\Omega_0)$ and $(J_1,\Omega_1)$ as in Proposition \ref{cpct} and generic paths $(J_t,\Omega_t)$ connecting them,
the space
\[\mathcal{PM}_{0}(\alpha_1,\ldots,\alpha_n)=\{(t,s)|s\in\mathcal{M}^{J_t,\Omega_t}_{0}(\alpha_1,\ldots,\alpha_n)\}\]
is a compact one-dimensional manifold. \end{prop}

\begin{proof} This follows immediately from the above discussion, noting that in the proof of Proposition \ref{cpct} we saw that any possible
boundary components of
$\mathcal{M}^{J}_{0}(\alpha_1,\ldots,\alpha_n)$ have real
codimension 2 and so will not appear in our one-dimensional
parametrized moduli space. \end{proof}

Note that we can orient these moduli spaces by using the spectral
flow of the linearization of the $\overline{\partial}$ operator at
an element
$s\in\mathcal{M}^{J,\Omega}_{0}(\alpha_1,\ldots,\alpha_n)$ acting
on sections of $s^{*}T^{vt}X_r (f)$ which preserve the incidence
conditions and the tangencies to $\Delta$;
$\mathcal{PM}_{0}(\alpha_1,\ldots,\alpha_n)$ will then be an
oriented cobordism between
$\mathcal{M}^{J_0,\Omega_0}_{0}(\alpha_1,\ldots,\alpha_n)$ and
$\mathcal{M}^{J_1,\Omega_1}_{0}(\alpha_1,\ldots,\alpha_n)$.
Accordingly, we may make the following definition.

\begin{dfn} Let $\alpha=\alpha_1 +\cdots +\alpha_n$ be a decomposition of $\alpha\in H^2(X,\mathbb{Z})$ which satisfies Assumption \ref{ass}.  Then \[
\widetilde{\mathcal{DS}_f}(\alpha;\alpha_1,\ldots,\alpha_n) \] is
defined as the number of points, counted with sign according to
orientation, in the space
$\mathcal{M}^{J,\Omega}_{0}(\alpha_1,\ldots,\alpha_n)$ for generic
$(J,\Omega)$ as in Proposition \ref{cpct}. \end{dfn}

\begin{thm} \label{tildesame} If $\alpha=\alpha_1+\cdots\alpha_n$ is a decomposition satisfying Assumption \ref{ass} then \[
\frac{(\sum d(\alpha_i))!}{\prod
(d(\alpha_i)!)}Gr(\alpha;\alpha_1,\ldots,\alpha_n)=\widetilde{\mathcal{DS}_f}(\alpha;\alpha_1,\ldots,\alpha_n)\]
provided that the degree of the fibration is large enough that
$\langle [\omega_{X}],[\Phi]\rangle>[\omega_{X}]\cdot \alpha$.
\end{thm}

\begin{proof}
Let $j$ be an almost complex structure on $X$ generic among those
compatible with the fibration $f\co X\to S^2$, and $\Omega$ a
generic set of $\sum d(\alpha_i)$ points.  The curves in $X$
contributing to $Gr(\alpha;\alpha_1,\cdots,\alpha_n)$ are unions
\[ C=\bigcup_{i=1}^{n} C^i \] of embedded $j$-holomorphic curves
$C^i$ which are Poincar\'e dual to $\alpha_i$ (note that
Assumption \ref{ass} implies that none of these curves will be
multiple covers) with $\Omega_i\subset C^i$ for some fixed generic
sets $\Omega_i$ of $d(\alpha_i)$ points. In Section 3 of
\cite{Usher} it was shown that there is no loss of generality in
assuming that $j$ is integrable near $\cup_i Crit(f|_{C^i})$, so
let us assume that this is the case. Where $s_C$ is the section of
$X_r (f)$ tautologically corresponding to $C$, in the context of
\cite{Usher} this local integrability condition was enough to
ensure that the almost complex structure $\mathbb{J}_j$ on
$X_r(f)$ constructed from $j$ was smooth on a neighborhood of
$s_C$.  Here that is not quite the case, for $\mathbb{J}_j$ might
only be H\"older continuous at the points of $Im(s_C)$
tautologically corresponding to the intersection points of the
various $C^i$.

However, just as in Section 5 of \cite{Usher}, we can still define
the contribution $r'(C)$ to $\widetilde{DS}_f
(\alpha_1,\ldots,\alpha_n)$ by perturbing $\mathbb{J}_j$ to a
generic almost complex structure $J$ which is compatible with the
strata and H\"older-close to $\mathbb{J}_j$, and then counting
with sign the elements of
$\mathcal{M}^{J,\Omega}_{0}(\alpha_1,\ldots,\alpha_n)$ which lie
near $s_C$; since the curves $C$ which contribute to
$Gr(\alpha_1,\ldots,\alpha_n)$ are isolated, and since the members
of
$\mathcal{M}^{\mathbb{J}_j,\Omega}_{0}(\alpha_1,\ldots,\alpha_n)$
are precisely the $s_C$ corresponding to the curves $C$, it
follows from Gromov compactness that for sufficiently small
perturbations $J$ of $\mathbb{J}_j$ all elements of
$\mathcal{M}^{J,\Omega}_{0}(\alpha_1,\ldots,\alpha_n)$ will be
close to one and only one of the $s_C$.  Thus \[
\widetilde{\mathcal{DS}}_f (\alpha_1,\ldots,\alpha_n)=\sum_{\pi\in
p(\Omega)}\sum_{C\in\mathcal{M}^{j,\Omega,\pi}(\alpha_1,\ldots,\alpha_n)}r'(C)
\] where $p(\Omega)$ is the set of partitions of $\Omega$ into
subsets $\Omega_i$ of cardinality $d(\alpha_i)$ and, writing
$\pi=(\Omega_1,\ldots,\Omega_n)$,
 $\mathcal{M}^{j,\Omega,\pi}(\alpha_1,\ldots,\alpha_n)$ is the
space of curves $C=\cup C^i$ contributing to
$Gr(\alpha;\alpha_1,\ldots,\alpha_n)$ with $C^i$ passing through
$\Omega_i$. Meanwhile, for any $\pi$, we have
\[ Gr(\alpha;\alpha_1,\ldots,\alpha_n)=\sum_{C\in
\mathcal{M}^{j,\Omega,\pi}(\alpha_1,\ldots,\alpha_n)} r(C),\]
$r(C)$ being the product of the spectral flows of the
linearizations of $\overline{\partial}_j$ at the embeddings of the
$C^i$ where $C=\cup C^i$.  The theorem will thus be proven if we
show that $r'(C)=r(C)$, which we now set about doing.

So let $C=\cup C^i\in
\mathcal{M}^{j,\Omega,\pi}(\alpha_1,\ldots,\alpha_n)$.  Taking $j$
generically, we may assume that all intersections of the $C^i$ are
transverse and occur away from $crit(f|_{C^i})$ (this follows from
the arguments of Lemma 2.1 of \cite{Usher}).  Let $p\in C^i\cap
C^k$.  In a coordinate neighborhood $U$ around $p$, where $w$ is a
holomorphic coordinate on the fibers and $z$ the pullback of the
coordinate on $S^2$, we may write \[ C^i\cap U=\{w=g(z)\} \qquad
C^k\cap U=\{w=h(z)\}.\]  If the almost complex structure $j$ is
given in $U$ by \begin{equation} \label{T01} T^{0,1}_j =\langle
\partial_{\bar{z}}+b(z,w)\partial_w,\partial_{\bar{w}}\rangle\end{equation}
(note that we may choose the horizontal tangent space so that
$b(0,0)=0$), that $C^i$ and $C^k$ are $j$-holomorphic amounts to
the statement that \[
\partial_{\bar{z}}g(z)=b(z,g(z)) \qquad \partial_{\bar{z}}h(z)=b(z,h(z));\] in particular, we have $g_{\bar{z}}(0)=h_{\bar{z}}(0)=0$.
Since  $C^i\pitchfork C^k$, we have $(g-h)_z (0)\neq 0$, and by
the inverse function theorem $(g-h)\co \mathbb{C}\to\mathbb{C}$ is
invertible on some disc $D_{2\delta}(0)$.  Let $g_t$ and $h_t$
$(t\in [0,1])$ be one-parameter families of functions satisfying
\begin{itemize} \item[(i)] $g_0=g$, $h_0=h$; \item[(ii)] On
$D_{2\delta}(0)$, $g_t-h_t$ is invertible as a complex-valued
smooth function, with inverse $p_t$; \item[(iii)] $g_t$ and $h_t$
agree with $g$ and $h$, respectively, outside $D_{2\delta}(0)$;
\item[(iv)] $g_t(0)=h_t(0)=\partial_{\bar{z}}g_t
(0)=\partial_{\bar{z}}h_t (0)=0$; and \item[(v)] $g_1(z)$ and
$h_1(z)$ are both holomorphic on $D_{\delta}(0)$.\end{itemize} Let
\[ C_{t}^{i}=(C^i\cap (X\setminus U))\cup\{w=g_t(z)\} \mbox{ and
} C_{t}^{k}=(C^k\cap (X\setminus U))\cup\{w=h_t(z)\}.\]  Now set
\[ b_t
(z,w)=(\partial_{\bar{z}}h_t)(z)+\partial_{\bar{z}}\left(g_t-h_t\right)\left(p_t(w-h_t(z))\right).\]

Then, since $p_t=(g_t-h_t)^{-1}$, \[ b_t
(z,h_t(z))=\partial_{\bar{z}}h_t(z)+\partial_{\bar{z}}(g_t-h_t)(0)=\partial_{\bar{z}}h_t
(z)\] while
\[ b_t (z,g_t(z))=\partial_{\bar{z}}h_t(z)+\partial_{\bar{z}}(g_t-h_t)(z)=\partial_{\bar{z}}g_t(z).\]

Let $b'_t$ agree with $b_t$ near $\{(z,w)\in C_{t}^i\cup
C_{t}^k|z\in D_{2\delta}(0)\}$ and with $b$ sufficiently far from
the origin in $U$.  Then defining $j'_t$ by
$T^{0,1}_{j'_t}=\langle
\partial_{\bar{z}}+b'_t\partial_w,\partial_{\bar{w}}\rangle$,
$j'_t$ agrees with $j$ near $\partial U$ and makes $C^{i}_{t}\cup
C^{k}_{t}$ holomorphic.  Further, we see that $b_1(z,w)\equiv 0$
for $z\in D_{\delta}(0)$, from which a Nijenhuis tensor
computation shows that $j'_1$ is integrable on a neighborhood of
the unique point $p$ of $C^{i}_{1}\cap C^{k}_{1}\cap U$.

Carrying out this construction near all intersection points of the
$C^i$, we obtain curves $C_t=\cup C^{i}_{t}$ and almost complex
structures $j'_t$ on $X$ such that $j'_1$ is integrable near all
intersection points of the $C^{i}_{1}$.  Since $j'_1$ agrees with
$j$ and $C^{i}_{1}$ with $C^i$ away from small neighborhoods of
these intersection points, $j'_1$ is also integrable on a
neighborhood of $crit(f|_{C^{1}_i})$ for each $i$.

If $p$ is a point of $C_1$ near which $j'_1$ is not already
integrable, then in a neighborhood $U$ of $p$ we have $C_1 \cap
U=\{w=g(z)\}$, and so the condition for an almost complex
structure $j'$ given by
$T^{0,1}_{j'}=\langle\partial_{\bar{z}}+b\partial_w,\partial_{\bar{w}}\rangle$
to make $C_1$ holomorphic near $p$ is just that
$\partial_{\bar{z}}g(z)=b(z,g(z))$, while the condition for $j'$
to be integrable in the neighborhood is that
$\partial_{\bar{w}}b(z,w)=0$. As in Lemmas 4.1 and 4.4 of
\cite{Usher}, then, we may easily find a path of almost complex
structures $j'_t$ ($1\leq t\leq 2$) such that each $j'_t$ makes
$C_1$ holomorphic and $j'_2$ is integrable on a neighborhood of
$C_1$.  So, changing notation slightly, we have proven:
\begin{lemma} \label{isotopy} There exists an isotopy $(C_t,j_t)$ of pairs consisting of almost complex structures $j_t$ compatible with the fibration $f\co X\to S^2$ and
$j_t$-holomorphic curves $C_t$ such that $(C_0,j_0)=(C,j)$ and
$j_1$ is integrable on a neighborhood of $C_1$. \end{lemma}

In the situation of the above lemma, $\mathbb{J}_{j_1}$ is not
only smooth but also integrable on a neighborhood of $C_1$; Lemma
4.2 of \cite{Usher} shows that if $j_1$ is chosen generically
among almost complex structures which make both $C_1$ and $f$
pseudoholomorphic and are integrable near $C_1$ the linearization
of $\bar{\partial}_{\mathbb{J}_{j_1}}$ at $s_C$ will be
surjective, as will the linearizations of $\bar{\partial}_{j_1}$
at the embeddings of each of the $C_{1}^{i}$. We now fix the
isotopy $C_t$ and the almost complex structure $j_1$ which is
nondegenerate in the above sense; Lemma \ref{isotopy} then gives a
path $j_t$ from $j=j_0$ to $j_1$ such that each $C_t$ is
$j_t$-holomorphic.  We may then define $r'_{j_t}(C_t)$ in the same
way as $r'(C)$, by counting $J$-holomorphic sections close to
$s_{C_t}$ for some $J$ H\"older-close to $\mathbb{J}_{j_t}$.
Meanwhile, if the linearization $D\bar{\partial}_{j_t}$ is
surjective at the embeddings of the $C^{i}_{t}$, its spectral flow
gives a number $r_{j_t}(C_t)$, and our goal is to show that
$r_{j_0}(C_0)=r'_{j_0}(C_0)$.  To this end, we see from Lemma 5.5,
Corollary 5.6, and their proofs in \cite{Usher} that:
\begin{lemma} For generic paths $j_t$ from $j_0$ to $j_1$ as above such that $C_t$ is $j_t$-holomorphic, the following statements hold.  $D\bar{\partial}_{j_t}$
is surjective at the embeddings of the $C^{i}_{t}$ for all but
finitely many values of $t$.  For $t$ near any value $t_0$ for
which  $D\bar{\partial}_{j_{t_0}}$ fails to be surjective, the set
of elements of
$\mathcal{M}^{j_t,\Omega}(\alpha_1,\ldots,\alpha_n)$ in a tubular
neighborhood of $C_t$ is given by $\{C_t,\tilde{C}_t\}$ for a
smooth family of curves $\tilde{C}_t$ with
$\tilde{C}_{t_0}=C_{t_0}$.  Further, for small $\epsilon >0$, we
have
\[ r'_{j_{t_0+\epsilon}}(C_{t_0+\epsilon})=r'_{j_{t_0-\epsilon}}(\tilde{C}_{t_0-\epsilon})=-r'_{j_{t_0-\epsilon}}(C_{t_0-\epsilon}) \] and \[
 r_{j_{t_0+\epsilon}}(C_{t_0+\epsilon})=r_{j_{t_0-\epsilon}}(\tilde{C}_{t_0-\epsilon})=-r_{j_{t_0-\epsilon}}(C_{t_0-\epsilon}). \]  Moreover, on intervals
not containing any $t_0$ for which $j_{t_0}$ has a non-surjective
linearization, $r'_{j_t}(C_t)$ and $r_{j_t}(C_t)$ both remain
constant.  \end{lemma}

Since (for generic paths $j_t$), $r'_{j_t}(C_t)$ and
$r_{j_t}(C_t)$ stay constant except for finitely many points at
which they both change sign, to show that
$r'_{j_0}(C_0)=r_{j_0}(C_0)$ it is enough to see that
$r'_{j_1}(C_1)=r_{j_1}(C_1)$.  But since $j_1$ is
\emph{integrable} and nondegenerate near $C_1$, as is
$\mathbb{J}_{j_1}$ near $s_{C_1}$, we immediately see that
$r'_{j_1}(C_1)=r_{j_1}(C_1)=1$, and the theorem follows.
\end{proof}

\begin{rem}  The above proof suggests a simplification of the
proof that $\mathcal{DS}=Gr$ in \cite{Usher}.  As mentioned above,
in Section 3 of \cite{Usher} it is shown that we can take the
almost complex structure $j$ to be integrable on neighborhoods of
the critical points of the various $f|_C$ for $C$ contributing to
$Gr(\alpha)$.  Given arbitrary generic fibration-compatible $j$,
however, as in the proof of Theorem \ref{tildesame}, the arguments
of Sections 4 and 5 of \cite{Usher} go through as long as we can
find an isotopy  $(C_t,j_t)$ of pairs consisting of almost complex
structures $j_t$ compatible with the fibration $f\co X\to S^2$ and
$j_t$-holomorphic curves $C_t$ such that $(C_0,j_0)=(C,j)$ and
$j_1$ is integrable on a neighborhood of $C_1$.  This is indeed
possible; if near a critical point of $f|_C$ $C$ has the form
$\{z=w^n+O(n+1)\}$, we can take $C_t$ such that $C_t$ agrees with
$C$ away from a neighborhood of $Crit(f|_C)$ and $C_1$ has the
form $\{z=w^n\}$ on a smaller neighborhood of the critical point,
and then we can choose $j_t$ to make $C_t$ holomorphic. (The
easiest approach to this seems to be to have $C_t$ be constant for
$t\leq 1/2$ and arrange the function $b_{1/2}(z,w)$ in the
notation (\ref{T01}) to depend only on $w$ near the critical
points; then for $t>1/2$, the form of $C_t$ determines uniquely a
$z$-independent function $b_t$ which causes $C_t$ to be
$j_t$-holomorphic, and we will have $b_1(z,w)=0$ near the critical
point.  Details are left to the reader.)
\end{rem}

\section{The family standard surface count} \label{secfam}

While much is known about the structure the Gromov--Taubes
invariants, which count embedded holomorphic curves in symplectic
4-manifolds, we know comparatively little about invariants
counting singular curves.  We explain here an approach to nodal
curves using Donaldson and Smith's constructions.

We should mention first of all that whereas Taubes' work gives us
a natural invariant $Gr(\alpha)$ counting all embedded curves
(regardless of their connected-component decomposition) Poincar\'e
dual to some class $\alpha$, if we instead wish to assemble all of
the possibly-reducible curves Poincar\'e dual to $\alpha$ and
having some number $n>0$ of ordinary double points into an
invariant $Gr_n(\alpha)$, it is somewhat unclear how we should
proceed in many cases.  Just as with the difficulties surrounding
the Gromov--Taubes invariant, this stems from the multiple-cover
problem: if for some class $\beta\in H^2(X,\mathbb{Z})$ and $m>1$
we have $d(\beta)\geq \max\{0, d(m\beta)-n\}$, then for generic
almost complex structures $j$ there will arise the possibility of
a sequence of curves Poincar\'e dual to $m\beta$ which have $n$
double points converging to an $m$-fold cover of a curve
Poincar\'e dual to $\beta$.  When $n=0$, as was noted in the
previous section the formula for $d(\beta)$ and the adjunction
formula imply that this only arises when $\beta$ is Poincar\'e
dual to a square-zero torus, and Taubes' work shows how to
incorporate multiple covers into the definition of $Gr$ in the
correct way.  When $n>0$, the equation $d(\beta)\geq d(m\beta)-n$
becomes easier to satisfy and it is less clear how multiple covers
should be dealt with, especially in the case of a strict
inequality $d(\beta)>d(m\beta)-n$, when the multiple covers form a
space of larger dimension than that of the space we are interested
in.

Of course, there will typically be at least some classes for which
this issue does not arise:

\begin{dfn} \label{sem} A
class $\alpha\in H^2(X,\mathbb{Z})$ is called \emph{strongly}
$n$-\emph{semisimple} if there exist \emph{no} decompositions
$\alpha= \alpha_1+\cdots+\alpha_l$ into nonnegatively-intersecting classes $\alpha_i$ such that each $\alpha_i$ has
$d(\alpha_i)\geq 0$ and is Poincar\'e dual to the image of a
symplectic immersion, and $\alpha_1$ is equal to $m\beta$ ($m>1$)
where $\beta$ satisfies $d(\beta)\geq
\max\{0,d(\alpha_1)-n+\alpha_1\cdot(\alpha-\alpha_1)\}$.  $\alpha$
is called \emph{weakly $n$-semisimple} if the only decompositions
$\alpha=\alpha_1+\cdots +\alpha_n$ as above which exist have
$\alpha_{1}^{2}=\kappa_X\cdot\alpha_1=0$.
\end{dfn}

For instance, every class is weakly $0$-semisimple, while the only
classes which are not weakly $1$-semisimple are those classes
$\alpha$ such that there exists a class $\beta\in
H^2(X;\mathbb{Z})$ such that $\beta\cdot(\alpha-2\beta)=0$ and
$\beta$ is Poincar\'e dual either to a symplectic sphere of square
0 or a symplectic genus-two curve of square 1, while
$\alpha-2\beta$ is Poincar\'e dual to some embedded (and possibly
disconnected) symplectic submanifold.  For strong semisimplicity,
one needs to add the assumption that $\alpha$ is not Poincar\'e
dual to a symplectic immersion having a component which is a
square-zero torus in a non-primitive homology class.

For a weakly- or strongly-$n$-semisimple classes $\alpha$, there
is an obvious analogue of the Gromov--Taubes invariant
$Gr_{n}(\alpha)$, defined by counting $j$-holomorphic curves $C$
which are unions of curves $C_i$ Poincar\'e classes $\alpha_i$
carrying multiplicities $m_i$ which are equal to 1 unless $C_i$ is
a square-zero torus with $\sum m_i\alpha_i=\alpha$, such that $C$
has $n$ transverse double points and passes through a generic set
of $d(\alpha)-n$ points of $X$; each such $C$ contributes the
product of the Taubes weights $r(C_i,m_i)$ to the count
$Gr_n(\alpha)$.  Since the condition of $n$-semisimplicity is
engineered to rule out the only additional possible source of
noncompactness of the relevant moduli spaces, the proof that
$Gr(\alpha)$ is independent of the choice of almost complex
structure used to define it goes through to show the same result
for $Gr_n(\alpha)$.

For that matter, if $\alpha$ is weakly $n$-semisimple and we have
$n_i\geq 0$ and $\alpha_i$ with $\sum \alpha_i=\alpha$ and $\sum
n_i=n-\sum_{i<j}\alpha_i\cdot\alpha_j$, we can form a refinement
$Gr_{(n_1,\ldots,n_k)}(\alpha;\alpha_1,\ldots,\alpha_k)$ along the
lines of Definition \ref{invts} which counts (modulo the usual
square-zero torus issues) curves with reducible components which
are Poincar\'e dual to the $\alpha_i$ and have $n_i$ transverse
self-intersections.  In this case, under Assumption \ref{ass} it
is also straightforward to modify the constructions of the
previous section to produce an invariant
$\widetilde{\mathcal{DS}}_{(n_1,\ldots,n_k)}(\alpha;\alpha_1,\ldots,\alpha_k)$
which counts holomorphic sections $s$ of $X_r(f)$ in the homotopy
class $c_{\alpha}$ which decompose into a sum of $C^1$ sections
$s_i\in c_{\alpha_i}$ such that each $s_i$ has $n_i$ tangencies to
the diagonal stratum of $X_{r_i}(f)$ and does not itself decompose
as a nontrivial sum of $C^1$ sections.  Furthermore, the proof of
Theorem \ref{tildesame} goes through unchanged to show that \[
Gr_{(n_1,\ldots,n_k)}(\alpha;\alpha_1,\ldots,\alpha_k)=\widetilde{\mathcal{DS}}_{(n_1,\ldots,n_k)}(\alpha;\alpha_1,\ldots,\alpha_k).\]

Instead, though, we aim to produce an invariant similar to
$Gr_{n}(\alpha)$ which does not require $\alpha$ to be
$n$-semisimple.  For general $\alpha$, the multiple cover problem
discussed above has its mirror on the side of
$\widetilde{\mathcal{DS}}$ in the fact that the moduli spaces for
the latter will tend to have undesirably-large components
consisting of sections which are mapped entirely into the diagonal
stratum, so $\widetilde{\mathcal{DS}}$ will not be much help
toward this goal.  Instead, we take a hint from the approach used
by A.K. Liu in \cite{Liu} and construct family versions of the
standard surface count.  These new invariants will use almost
complex structures which generally do not make the diagonal
stratum pseudoholomorphic, and so we will not encounter moduli
spaces with unexpectedly large components consisting of sections
mapped into $\Delta$.

Be given a symplectic Lefschetz fibration $f\co X\to S^2$. Write
$f_0=f$, $X_0=\{pt\}$, $X_1=X$, and let $g_0\co X_1\to X_0$ be the
map of $X$ to a point.  As in \cite{Liu}, for $n\geq 1$ form
$X^{0}_{n+1}=X_{n}\times_{g_{n-1}}X_{n}$, and let $X_{n+1}$ be the
blowup of the relative diagonal in $X_{n+1}^{0}$. Let $g_{n}\co
X_{n+1}\to X_{n}$ be the projection onto the first factor.  Each
$X^b:=g_{n}^{-1}(b)$ ($b\in X_{n}$)  is then an $n$-fold blowup of
$X$, with the parameter $b$ indicating which points have been
blown up. Composing the maps $g_n$ gives a map $X_{n+1}\to X_1=X$;
let $f_n\co X_{n+1}\to S^2$ be the composition of this map with
the Lefschetz fibration $f$. (Equivalently, on each $n$-fold
blowup $X^b=g_{n}^{-1}(b)$, $f_n|_{X^b}$ is the composition of the
blowdown map with the Lefschetz fibration $f$.)

Write $f^{b}=f_{n}|_{X^b}$.  $f^b\co
X\#n\overline{\mathbb{C}P^2}\to S^2$ then has the same structure
as $f$, except that if $k$ points on some fiber (in class
$[\Phi]$) are among the blown up points, that (initially
irreducible) fiber has been replaced by a reducible curve with
components in classes $[\Phi]-E_1-\cdots -E_k$, $E_1,\ldots,E_k$,
where the $E_i$ are classes of exceptional spheres.
Straightforward local coordinate calculations show that, if none
of the blown-up points are critical points of any of the $f_i$
($i<n$), then the only intersection points between components are
ordinary double points, and that near the double points $f^b$ has
form $(z,w)\mapsto zw$. In particular, each $f^b=f_n|_{X^b}$ is
still a Lefschetz fibration provided that no critical points of
any of the intermediate fibrations are blown up in forming $X^b$.

\begin{notn} Denote a point $b\in X_n$ by $(p_1,\ldots,p_n)$,
where each $p_{j+1}\in X^{(p_1,\ldots,p_{j})}$.  Let:
\begin{enumerate}
\item[(i)] $X'_{n}$ be the set of $(p_1,\ldots,p_n)\in X_n$ such
that no $p_{j+1}$ is a critical point of
$f^{(p_1,\ldots,p_{j})}\co X^{(p_1,\ldots,p_j)}\to S^2$.
\item[(ii)] $X''_{n}$ be the set of $(p_1,\ldots,p_n)\in X_n$ such
that no $p_{j+1}$ lies in a singular fiber of
$f^{(p_1,\ldots,p_{j})}\co X^{(p_1,\ldots,p_j)}\to S^2$.
\end{enumerate}
\end{notn}

If $b\in X'_n$, then, our above remarks show that $f^b\co X^b\to
S^2$ is a Lefschetz fibration; if moreover $b\in X''_n$, then no
fiber of $f^b$ will contain more than one critical point (and also
none of the $n$ blowups involved in the creation of $X^b$ will be
at a point on an exceptional divisor of a previous blowup).

\begin{notn} \emph{(i)} For any $b\in X'_n$, $F^b\co X^{b}_{r}(f^b)\to S^2$ shall denote
the relative Hilbert scheme constructed from $f^b$ as in the
Appendix of \cite{DS} and Section 3 of \cite{Smith}.\newline
\emph{(ii)} $\mathcal{X}^{n}_{r}(f)=\{(D,b): b\in X'_n,\,D\in
X^{b}_{r}(f^b)\}$.  In particular we have a map $\mathcal{F}^n\co
\mathcal{X}^{n}_{r}(f)\to S^2\times X'_n$.
\end{notn}

For $b\in X''_n$, $X^b$ contains disjoint exceptional divisors
$E_1,\ldots,E_n$, and our intention is to define an invariant
counting sections of the various $X^{b}_{r}(f^b)$ which descend to
curves Poincar\'e dual to $\alpha-2\sum PD(E_i)$, as $b$ ranges
over $X''_n$.  We have to be somewhat careful in the definition of
this invariant, though, since our parameter space $X''_n$ is
noncompact.

\begin{lemma} For $b\in X''_n$, $X^{b}_{r}(f^b)$ is a smooth
symplectic manifold, as is the total space of
$\mathcal{X}^{n}_{r}(f)\to S^2\times X'_n$.\end{lemma}
\begin{proof} That the relative Hilbert scheme constructed from
any Lefschetz fibration (such as $f^b$ when $b\in X''_n$) in which
there is at most one critical point per fiber is smooth is shown
in Theorem 3.4 of \cite{Smith} (as noted in Remark 3.5 of
\cite{Smith}, Smith's provision of a local coordinate description
for the relative Hilbert scheme makes irrelevant his assumption
that all of the fibers of the Lefschetz fibration are
irreducible).  When $b\in X'_{n}\setminus X''_n$, so that $f^b$,
while still a Lefschetz fibration, may have more than one critical
point per fiber, the individual $X^{b}_{r}(f^b)$ will tend not to
be smooth near points on the Hilbert scheme of the singular fibers
$\Sigma_0$ which are sent by the map $Hilb^{[r]}\Sigma_0\to
S^r\Sigma_0$ to divisors which contain more than one of the nodes
of $\Sigma_0$.  We will show presently, though, that the freedom
to vary $b\in X'_{n}$ results in the total space
$\mathcal{X}^{n}_{r}(f)$ still being smooth at these points.

To see this, note that Donaldson and Smith show (c.f. the proof of
Proposition A.8 of \cite{DS}) that when $f$ only has one node per
fiber, at a singular point of a fiber of $X_s(f)$ (corresponding
to a divisor with points near the node of a fiber) the behavior of
$F\co X_s(f)\to S^2$ is modeled by $(z_1,\ldots,z_{s+1})\mapsto
z_1z_2$.  When there are multiple nodes in a fiber, then, the
relative Hilbert scheme will be modeled near a point corresponding
to a divisor containing $s_i$ copies of the nodes $p_i$
($i=1,\ldots,l$) by the fiber product of the various maps
$(z^{(i)}_1,\ldots,z^{(i)}_{s_i+1})\mapsto
z^{(i)}_{1}z^{(i)}_{2}$.  This fiber product is the common
vanishing locus of the various
$z^{(i)}_{1}z^{(i)}_{2}-z^{(j)}_{1}z^{(j)}_{2}$ (which is of
course singular where $z^{(i)}_{1}=z^{(i)}_{2}=0$ for all $i$).

More generally, though, if $p_i$ is a node lying near the fiber
over zero, $X_{s}(f)\to S^2$ is modeled near points corresponding
to divisors with points near $p_i$ by \\
$(z^{(i)}_1,\ldots,z^{(i)}_{s})\mapsto
z^{(i)}_{1}z^{(i)}_{2}+f(p_i)$.  In our present context the
fibration map is $f^b$; say for notational simplicity that
$b=(p_1,\ldots,p_n)$ gives rise to an $n$-fold blowup with all
exceptional divisors in the same fiber (of course if some
exceptional divisors are in different fibers we can work
fiber-by-fiber and reduce to this case). The space
$\mathcal{X}^{n}_{r}(f)$ is then, at worst, modeled locally by
\begin{equation} \label{model}
\{(\vec{z}^{(0)},\vec{z}^{(1)},\ldots,\vec{z}^{(n)},q_1,\ldots,q_n):
z^{(0)}_{1}z^{(0)}_{2}=z^{(i)}_{1}z^{(i)}_{2}+f^{(p_1,\ldots,p_{i-1})}(q_i)\}.
\end{equation}  Here
 $\vec{z}^{(0)}$ are the coordinates on the relative Hilbert scheme corresponding to divisors which contain any nodes that may have existed in our fiber before blowing up (and we are of course assuming throughout that the original $f$ was chosen so that there is at most one such). The
$q_i$ are elements of a coordinate chart centered on $p_i\in
X^{(p_1,\ldots,p_{i-1})}$.  But (\ref{model}) defines a smooth
manifold at any point with $q_i=p_i$ as long as none of the $p_i$
are  critical points for $f^{(p_1,\ldots,p_{i-1})}$, and this
latter condition is precisely ensured by the fact that $b\in
X'_n$.

This shows that $\mathcal{X}^{n}_{r}(f)$ is smooth; the existence
of a symplectic structure on it then follows exactly as in the
proof of the existence of a symplectic structure on $X_r(f)$ in
\cite{DS}: where $\mathcal{X}^{n}_{r}(f)$ fails to be a fibration
we have a local K\"ahler model for it, and we can extend the
resulting form to the entire manifold by the usual methods of
Gompf and Thurston. \end{proof}

We consider almost complex structures $J$ on the $X^{b}_{r}(f^b)$
which make the fibration maps $F^b\co X^{b}_{r}(f^b)\to S^2$
pseudoholomorphic and have the following special type: for each
reducible fiber of $X^b$, letting $E$ denote the union of the
spherical components of that fiber, we require that there exist
neighborhoods $U\supset V$ of $E$ with $f^b(U)=f^b(V)=W\subset
S^2$ and almost complex structures $J_{1}^{q}$ and $J_{2}^{q}$ on
the restricted relative Hilbert schemes $X_q(f^b|_U)$ and
$X_{r-q}(f^b|_{(f^b)^{-1}(W)-V})$ such that the natural ``addition
map'' $X^{b}_{q}(f^b|_U)\times_{F^b}
X^{b}_{r-q}(f^b|_{(f^b)^{-1}(W)-V})\to X^{b}_{r}(f^b)$ is
$(J_{1}^{q}\times_{F^b}J_{2}^{q}, J)$-holomorphic; moreover, we
require that $J_{1}^{q}$ agree with the complex structure induced
(via the algebro-geometric description for the relative Hilbert
scheme given in Section 3 of \cite{Smith}) by an integrable
complex structure on $U\supset E$ with respect to which $f^b$ is
holomorphic.  Note that one way of forming such a $J$ is by taking
any almost complex structure on $X^{b}_{r}(f^b)$ which agrees near
the singular fibers with the almost complex structure
$\mathbb{J}_j$ tautologically corresponding to a structure $j$ on
$X^{b}$ which is integrable near the singular fibers of $X^b$.  If
$j$ is instead integrable only on the neighborhood $U$ of the
exceptional spheres, we still obtain a H\"older almost complex
structure satisfying the requirement, which may then be
H\"older-approximated by smooth almost complex structures also
satisfying the requirement by smoothing the almost complex
structures $J_{2}^{q}$ in a coherent way at points of the
$X_{r-q}(f^b|_{f^{-1}(W)-V})$ corresponding to divisors having
points missing $U$.

Let $\mathcal{J}$ denote the space of smooth tame almost complex
structures on $\mathcal{X}^{n}_{r}(f)$ which restrict to each
$X^{b}_{r}(f^b)=(\mathcal{F}^n)^{-1}(S^2\times\{b\})$ as a $J$ of
the above form. For each $b$, the blowdown map $\pi^b\co X^b\to X$
naturally induces a generically injective map $\Pi^b\co
X^{b}_{r}(f^b)\to X_r(f)$ on relative Hilbert schemes.  For $J\in
\mathcal{J}$ we obtain commutative diagrams
\[
\begin{CD}
{X^{b}_{q}(f^b|_U)\times_{F^b}
X^{b}_{r-q}(f^b|_{(f^b)^{-1}(W)-V})}@>>>{X^{b}_{r}(f^b)}\\
@VVV @ VV{\Pi^b}V\\
{X_{q}(f|_{\pi^b(U)})\times_{F}
X_{r-q}(f^b|_{f^{-1}(W)-\pi^b(V)})}@>>>{X_r(f)}\\
\end{CD}\]
in which $\Pi^b$ pushes $J$ forward to a smooth almost complex
structure $J_b$ on $X_r(f)$.  The $J_b$ vary smoothly in $b$, and
indeed extend by continuity to a smoothly $X_n$-parametrized
family of almost complex structures on $X_{r}(f)$ (rather than
just an $X'_n$-parametrized family).  Since our sections of the
$F^b\co X^{b}_{r}(f^b)\to S^2$ pass through all of the fibers of
$F^b$, restricting our almost complex structures to behave in this
way near the special fibers of $F^b$ will not prevent moduli
spaces of $J$-holomorphic sections of the $X^{b}_{r}(f^b)$ from
being of the expected dimension for generic $J\in \mathcal{J}$.

For $\alpha\in H^2 (X;\mathbb{Z})$, $b\in X''_n$, and $e_i$
($i=1,\ldots,n$) the Poincar\'e duals to the exceptional divisors
of the blowups which form $X^{b}$, note that the expected complex
dimension of the space of curves Poincar\'e dual to $\alpha-2\sum
e_i$ is $d(\alpha-2\sum e_i)=d(\alpha)-3n$, so since the the
\emph{real} dimension of $X''_n$ is $4n$ we would expect the space
of such curves appearing in any $X^b$ as $b$ ranges over $X''_n$
to have complex dimension $d(\alpha)-n$.

\begin{lemma} \label{cpct2} Let $\alpha\in H^2(X;\mathbb{Z})$, and choose a generic set $\Omega$ of $d(\alpha)-n$ points in $X$. For generic $J\in \mathcal{J}$, and also for generic paths $J_t$ in $\mathcal{J}$ connecting two such generic $J$, the spaces \[
\mathcal{M}^{n}_{J,\Omega}(\alpha-2\sum e_i)=\{(s,b): b\in
X''_n,\, s\in c_{\alpha-2\sum e_i}\subset
\Gamma(X^{b}_{r}(f^b)),\, \overline{\partial}_J
s=0,\,\Omega\subset C_s\} \] and \[
\mathcal{PM}^{n}_{(J_t),\Omega}(\alpha-2\sum e_i)=\{(s,b,t): b\in
X''_n,\, s\in c_{\alpha-2\sum e_i}\subset
\Gamma(X^{b}_{r}(f^b)),\, \overline{\partial}_{J_t}
s=0,\,\Omega\subset C_s\} \] are compact manifolds of real
dimensions zero and one, respectively, provided that $r=\langle
\alpha,[\Phi]\rangle \geq g+3n$ where $g$ is the genus of the
generic fiber of $f\co X\to S^2$.
\end{lemma}

\begin{proof} That the dimensions will generically be as expected is a standard result (for the general theory of ``parametrized Gromov--Witten invariants''
of the sort that we are in the process of defining see
\cite{Ruan}, though the compactness result proved presently makes
much of Ruan's machinery unnecessary for our purposes), so we only
concern ourselves with compactness.

Let $(s^m,b^m)$ be a sequence of $J$-holomorphic sections (or
$J_{t_m}$-holomorphic sections with $J_{t_m}\to J$) from either of
the sets at issue.  \emph{A priori}, there are two possible
sources of noncompactness: the $b^m$ might have a limit in
$X_n\setminus X''_n$, or the $b^m$ might converge to $b\in X''_n$
with the $s^m$ converging to a bubble tree.  As usual for
section-counting invariants, we can eliminate the second
possibility: because $J|_{X^{b}_{r}(f)}$ makes  $X^{b}_{r}(f)\to
S^2$ holomorphic, any bubbles must be contained in the fibers, and
so the section component of the resulting bubble tree would
descend to a set Poincar\'e dual to $\alpha-2\sum e_i-PD(i_*B)$,
where $B$ is some class in one of the fibers $(f^b)^{-1}(t)$ of
the fibration $f^b:X^b\to S^2$.  If $(f^b)^{-1}(t)$ is
irreducible, $B$ will necessarily be a positive multiple of the
fundamental class of the fiber, and just as in Section 4 of
\cite{Smith} we will have $d(\alpha-2\sum e_i-PD(i_*B))\leq
d(\alpha-2\sum e_i)-(r-g+1)$, which rules such bubble trees out
for generic one-parameter families of $J$.
 If $(f^b)^{-1}(t)$ is reducible, with components in classes $[\Phi]-E$ and $E$, then $B$ will have form $m([\Phi]-E)+pE$ where $m,p\geq 0$ and at least one is positive,
and a routine computation then yields that \[ d(\alpha-2\sum
e_i-PD(i_*B))-d(\alpha-2\sum
e_i)=-m(r-g+1)-\frac{5}{2}(p-m)-\frac{1}{2}(p-m)^2,\] which, since
we have assumed that $r\geq g+3$, will always be negative when
$m,p\geq 0$ and are not both zero. Thus for generic $J$ or $J_t$,
none of the possible bubble trees appear.

There remains the issue that the $b^m$ might converge to some
$b\notin X''_n$.  We rule this out in two steps: first, we prove:
\begin{sub} If $b\in X_n\setminus X'_n$ then $b^m$ cannot converge
to $b$. \end{sub} \begin{proof}[Proof of the sublemma] Let
$\pi^{b^m}\co X^{b^m}\to X$ be the blowdown map, and
\\ $\Pi^{b^m}\co X^{b^m}_{r}(f^{b^m})\to X_r(f)$ the map that it induces
on relative Hilbert schemes.  By the definition of our space
$\mathcal{J}$ of almost complex structures, the $\Pi^{b^m}\circ
s^m$ are $J_{b^m}$-holomorphic sections of $X_r(f)$ in the class
$c_{\alpha}$, and so converge modulo bubbling to a
$J^b$-holomorphic section $\bar{s}$ of $X_r(f)$.  In fact, we can
rule out bubbling, since we can assume that the family $J_b$ is
regular as a $4n$-real-dimensional family of almost complex
structures on $X_r(f)$, and so as above no bubbles can form in the
limit thanks to the fact that all fibers of $f$ are irreducible
and
\begin{align*}2n+d(\alpha-mPD[\Phi])&=d(\alpha)+2n-m(r-g+1)\\&\leq
d(\alpha)-n-(r-g+1-3n)<d(\alpha)-n \end{align*} by the hypothesis
of the lemma.

 Since $b\notin X'_n$,
where $b=(p_1,\ldots,p_n)$ there will be some minimal $l$ such
that $p_{l+1}$ is a critical point of $f^{(p_1,\ldots,p_{l})}\co
X^{(p_1,\ldots,p_l)}\to S^2$.
 Suppose first that $p_{l+1}$ lies on just one
irreducible component of its fiber (so that it is a double point
of that component).  Write
$t^m=f^{(p^{m}_{1},\ldots,p^{m}_{l})}(p^{m}_{l+1})$ and
$T=f^{(p_1,\ldots,p_l)}(p_{l+1})$. Now since $C_{s^m}\subset X^b$
meets the exceptional divisor formed by blowing up $p^{m}_{l+1}$
transversely exactly twice, we deduce that $\Pi\circ s^m\in
\Gamma(X_r(f))$ acquires a tangency to the diagonal at a divisor
containing two copies of $\pi^{b^m}(p_{l+1}^{m})$; more
specifically, assuming that $\bar{s}(T)$ corresponds to a divisor
containing $p_{l+1}$ with multiplicity $q$, for large $m$ in a
neighborhood $U$ around $T,t^m\in S^2$ we have a decomposition
$\Pi\circ s^m|_U=+(s^{m}_{1},s^{m}_{2})$ into disjoint summands
$s^{m}_{1}\in \Gamma(X_q(f)|_U)$ and $s^{m}_{2}\in
\Gamma(X_{r-q}(f)|_U)$, with $s^{m}_{1}$ tangent to the diagonal
at a point of form $\{p^{m}_{l+1},p^{m}_{l+1},x_3,\ldots,x_q\}$.
Since the divisors $s^{m}_{1}(t)$ and $s^{m}_{2}(t)$ are disjoint
for $t\in U$, the smoothness of the $\Pi\circ s^m$ implies the
smoothness of $s^{m}_{1}$ and $s^{m}_{2}$ over $U$.  Similarly,
where $V$ is a neighborhood of $p_{l+1}$ with $f(V)\subset U$
$\bar{s}$ splits near $T$ into disjoint sections $\bar{s}_1$ of
$\mathcal{H}_q\cong X_{q}(f|_V)$ and $\bar{s}_2$ of
$X_{r-q}(f|_{f^{-1}(f(V))-V})$; here $\mathcal{H}_q$ is the
$q$-fold relative Hilbert scheme of the map $(z,w)\mapsto zw$.
Moreover, we have $s^{m}_{1}\to \bar{s}_1$.  But then since
$p_{l+1}^{m}\to p_{l+1}$, $\bar{s}_1$ must then be tangent to the
diagonal in $\mathcal{H}_q$ at a point corresponding to
$\{(0,0),\ldots,(0,0)\}\in Sym^{q}\{zw=0\}$.  This, however, is
impossible, since $\bar{s}_1$ is a \emph{section} of
$\mathcal{H}_q$, so that $Im(d\bar{s}_1)_T$ cannot be tangent to
the fiber, whereas according to Theorem \ref{appmain} in Section
\ref{app1} the tangent cone to $\Delta\subset\mathcal{H}_q$ is
contained in the tangent space to the fiber at $\bar{s}_1(T)$.

The other possibility is that $p_{l+1}$ is an intersection point
between two irreducible components of its fiber, in which case one
of those components is the exceptional sphere $E$ formed by a
previous blowup (say at $p_a$).  Where again
$t^m=f^{(p^{m}_{1},\ldots,p^{m}_{l})}(p^{m}_{l+1})$, in local
coordinate systems $U^m$   around $t^m$ (which may be shrinking
but are scale-invariant) we have
\[ \Pi\circ s^m=\{c_m z,d_m z\}+s^{m}_{2}(z) \] where $s^{m}_{2}$
is a local section of $X_{r-2}(f)$ which does not meet
$z\mapsto\{c_m z,d_m z\}$.  Now the fact that $p^{l+1}_{m}\to
p_{l+1}$  which is an intersection point between the fiber
containing $p_{l+1}\in X^{(p_1,\ldots,p_l)}$ and the exceptional
sphere of one of the blowups implies that, in $X$ (where the
blowup has not yet taken place), the two branches $c_mz$ and
$d_mz$ of $\Pi\circ s^m$ near $\pi^{b^m}(p^{l+1})$ both tend
toward the vertical, so that $c_m,d_m\to \infty$.  But then this
implies that $|d(\Pi\circ s^m)_{t^m}|\to \infty$, which is
impossible by elliptic regularity since $\Pi\circ s^m\to \bar{s}$.
\end{proof}

Finally we show that, generically, if $b^m\to b\in X'_n$ then in
fact $b\in X''_n$.  Indeed, since $b\in X'_n$, so that
$X^{b}_{r}(f^b)\subset \mathcal{X}^{n}_{r}(f)$, Gromov compactness
on the symplectic manifold $\mathcal{X}^{n}_{r}(f)$ implies that
after passing to a subsequence the sections $s^m$ will converge to
some smooth section $\bar{s}$ of $X^{b}_{r}(f^b)$. Just as above,
the fact that $\bar{s}$ is a smooth section implies that it misses
the critical locus of $F^b\co X^{b}_{r}(f^b)\to S^2$; in
particular, if $b\in X'_n\setminus X''_n$, $Im(\bar{s})$ is
contained in the smooth part of the relative Hilbert scheme
$X^{b}_{r}(f^b)$.  But then a neighborhood of $Im(\bar{s})$ in
$X^{b}_{r}(f^b)$ will be diffeomorphic to a neighborhood of
$Im(s^m)$ in $X^{b^m}_{r}(f^{b_m})$ for large m, and so the index
of the Cauchy-Riemann operator acting on perturbations of the
former will be the same as the index of the Cauchy-Riemann
operator acting on perturbations of the latter, namely
$d(\alpha)-3n$. Hence since the real dimension of $X'_n\setminus
X''_n$ is $4n-2$, the expected complex dimension of the space of
possible limits $\bar{s}$ with $b\in X'_n\setminus X''_n$ is
$d(\alpha)-n-1$, so for generic $J$, and also for generic
one-real-parameter families $J_t$, on $\mathcal{X}^{n}_{r}(f)$, no
such limits $\bar{s}$ with $C_{\bar{s}}$ satisfying our
$d(\alpha)-n$ incidence conditions will exist. \end{proof}

Given this compactness result, the standard cobordism argument
permits us to make the following definition.
\begin{dfn} \label{fds} Let $\alpha$ be as in Lemma \ref{cpct2}.
$\mathcal{FDS}^{n}_{f}(\alpha-2\sum e_i)$ is then defined as the
number of elements, counted with sign according to the spectral
flow, in the moduli space $\mathcal{M}^{n}_{J,\Omega}(\alpha-2\sum
e_i)$ for generic $J$ and $\Omega$ as in Lemma \ref{cpct2}.
\end{dfn}

\begin{thm}\label{famsame} Suppose that $\alpha$ is as in Lemma \ref{cpct2} and is strongly $n$-semisimple.
Then
\[ n!Gr_n(\alpha)=\mathcal{FDS}^{n}_{f}\left(\alpha-2\sum
e_i\right),
\] provided that $\langle \omega_{X},[\Phi]\rangle>\omega_{X}\cdot \alpha\geq g(\Phi)+3n$.
\end{thm}
\begin{proof}
As in the proof of Theorem \ref{tildesame}, we may evaluate
$Gr_n(\alpha)$ using an almost complex structure $j$ which makes
the Lefschetz fibration $f$ pseudoholomorphic and which has the
property that, for any of the curves $C=\bigcup_i C^i$ being
counted by $Gr_n(\alpha)$, $j$ is integrable on a neighborhood of
$\bigcup_i Crit(f|_{C^i})$; each intersection point between the
$C^i$ occurs away from $\bigcup_i Crit(f|_{C^i})$; and $C$ misses
the critical locus of the fibration $f$. For each $b\in X_n$, let
$j_b$ be pullback of $j$ via the blowup $\pi^b\co X^b\to X$ (see
Section \ref{app2} for the proof that $j_b$ exists and is
Lipschitz), so that $X^b\to X$ is $(j_b,j)$-holomorphic. Then, for
any of the $n!$ elements $b$ of $X'_n$ corresponding to the $n!$
different orders in which the nodes of $C$ may be blown up, the
proper transform $\tilde{C}$ of $C$ will be a curve in  $X^b$
(with $b\in X'_n$ as a result of the fact that $C$ misses the
critical points of $f$) Poincar\'e dual to $\alpha-2\sum e_i$. In
fact, we claim that for a generic initial choice of $j$ these
proper transforms $\tilde{C}$ are guaranteed to be the only
$j_b$-holomorphic curves Poincar\'e dual to $\alpha-2\sum e_i$  in
any $X^b$ which have no components contained in the fibers of
$f^b\co X^b\to S^2$.

Indeed, suppose that $\tilde{C}=\cup_i\tilde{C}_i$ is a
$j_b$-holomorphic curve in one of the $X^b$ Poincar\'e dual to
$\alpha-2\sum e_i$, with the (possibly-multiply-covered)
components $\tilde{C}_i$ Poincar\'e dual to $\beta_i-\sum
c_{ik}e_k$.  We need to show that, where $\pi^b\co X^b\to X$ is
the blowup, $\pi^b(\tilde{C})$ has $n$ nodes, located at the
points $p_i,\ldots,p_n$ which were blown up to form $X^b$ (as
$\pi^b(\tilde{C})$ is obviously a $j$-holomorphic curve Poincar\'e
dual to $\alpha$). Now for each $k$, $\sum_i c_{ik}=-2$, while by
positivity of intersections in $X^b$, we have each $c_{ik}\leq 0$.
If $k$ is such that there are distinct $q$ and $s$ with
$c_{qk}=c_{sk}=-1$, then the curves $\pi(\tilde{C}_q)$ and
$\pi(\tilde{C}_s)$ intersect transversely at the point $p_k$,
contributing the desired node.  On the other hand, if $k$ is such
that the only nonzero $c_{ik}$ is some $c_{qk}=-2$, then
$\pi^b(C_q)$ might \emph{a priori} be either a singly-covered
curve Poincar\'e dual to $\beta_q$ which has a self-intersection
at $p_k$, or a double cover of a curve in class $\beta_q/2$ which
passes through $p_k$.  However, the $n$-semisimplicity condition
rules the second possibility out for generic choices of $j$, since
we will have either $d(\beta_q/2)<0$ or
$d(\beta_q/2)<d(\beta_q)-n\leq d(\alpha)-n$, and so no such curves
satisfying our incidence conditions will exist.

We conclude, then, that the only $j_b$-holomorphic curves
$\tilde{C}$ in any $X^b$ Poincar\'e dual to $\alpha-2\sum e_i$ are
proper transforms of $j$-holomorphic curves which contribute to
$Gr_n(\alpha)$.  With this established, the proof of the theorem
becomes almost just an application of our usual methods.  Since
the restriction of $j_b$ to the exceptional spheres is standard,
we can choose smooth almost complex structures $j'_b$ which are
integrable near the exceptional spheres and are $C^0$-close to the
$j^b$.  By Gromov compactness for $C^0$ convergence of almost
complex structures \cite{IS} and the fact that $d(\alpha-2\sum
e_i)=-n$, we deduce as usual that for generic choices of these
perturbed $j'_b$ each $\tilde{C}$ will have finitely many
$j'_{b_i}$-holomorphic curves $\tilde{C}_1,\ldots,\tilde{C}_N$
near it (for various $b_i$ near $b$). On the relative Hilbert
schemes we have almost complex structures $\mathbb{J}_{j'_b}$.  If
$\tilde{C}_i$ is one of the curves above with the intersections of
its components resolved by the blowup $X^{b_i}\to X$, we define
$r''(\tilde{C}_i)$   as the signed count of $J_{b'}$ holomorphic
sections of $X^{b'}_{r}(f^{b'})$ near $s_{\tilde{C}_i}$ for $b'$
near $b_i$ and $J_{b'}$ a generic family of smooth almost complex
structures H\"older-close to the $\mathbb{J}_{j_{b_i}}$.

For $C$ a curve contributing to the Gromov invariant with nodes
resolved by $X^b\to X$ and proper transform $\tilde{C}$, we define
the contribution $r'(C)$ of $C$ to $\mathcal{FDS}$ as
$\sum_{i=1}^{n} r''(\tilde{C}_i)$ where the $\tilde{C}_i$ are
obtained as above. When $j$ is integrable near $C$, each $j_{b'}$
will be integrable near $\tilde{C}$ and near the exceptional
spheres of $X^{b'}$ for $b'$ near $b$,, so that the first
perturbation of the $j_{b'}$ to $j'_{b'}$ is not necessary and the
only $\tilde{C}_i$ is $\tilde{C}$ itself.  Moreover, each
$\mathbb{J}_{j_{b'}}$ will be integrable near $s_{\tilde{C}}$ for
$b'$ near $b$, and so (under suitable nondegeneracy assumptions)
both contributions will be 1. Further, exactly as in the proof of
Theorem \ref{tildesame}, the contributions transform under
variations in $j$ in the same way by virtue of the fact that
$\mathcal{FDS}$ is independent of the almost complex structure
used to define it. The agreement of the invariants then follows.
\end{proof}

If $\alpha$ is only weakly $n$-semisimple, then if $C\in
PD(\alpha)$ is the disjoint union of a double cover of a
square-zero torus with a curve having $n-1$ nodes, then the proper
transform of $C$ under blowup at the nodes of $C$ and at any point
on the torus gives a curve in some $X^b$ Poincar\'e dual to
$\alpha-2\sum e_i$, even though $C$ does not contribute to
$Gr_n(\alpha)$. On perturbing the family $(\mathbb{J}_{j_b})$ on
$\mathcal{X}^{n}_{r}(f)$ to a generic family $(J_b)$, we might
find that the sections corresponding to these curves contribute to
$\mathcal{FDS}^{n}_{f}(\alpha-2\sum e_i)$.  It seems reasonable,
though, to believe that these additional contributions could be
expressed in terms of the various other Gromov invariants of $X$,
consistently with Conjecture \ref{conj}.

\section{A review of Smith's constructions} \label{review}

Our vanishing theorem for $\mathcal{FDS}$ will follow by adapting
the constructions found in Section 6 of \cite{Smith} to the family
context. Let us review these.

In addition to the relative Hilbert scheme, Donaldson and Smith
constructed from the Lefschetz fibration $f\co X\to S^2$ a
\emph{relative Picard scheme} $P_r(f)$ whose fiber over a regular
value $t\in S^2$ is naturally identified with the Picard variety
$Pic^r\Sigma_t$ of degree-$r$ line bundles on $\Sigma_t$.  Over
each $\Sigma_t$, we have an Abel--Jacobi map $S^r\Sigma_t\to
Pic^r\Sigma_t$ mapping a divisor $D$ to its associated line bundle
$\mathcal{O}(D)$; letting $t$ vary over $S^2$, we then get a map
\[ AJ\co X_r(f)\to P_r(f)\] (that all of these constructions
extend smoothly over the critical values of $f\co X\to S^2$ is
seen in the Appendix of \cite{DS}).  Meanwhile, by composing the
Abel--Jacobi map for effective divisors of degree $2g-2-r$ with
the Serre duality map $\mathcal{L}\mapsto\kappa_{\Sigma_t} \otimes
\mathcal{L}^{\vee}$, we obtain a map \begin{align} i\co
X_{2g-2-r}(f)&\to P_r (f) \nonumber\\ D&\mapsto
\mathcal{O}(\kappa-D). \end{align} Moreover, using a result from
Brill-Noether theory due to Eisenbud and Harris \cite{EH}, Smith
obtains that (cf. Theorem 6.1 and Proposition 6.2 of
\cite{Smith}):
\begin{lemma}(\cite{Smith})\label{BN}
For a generic choice of fiberwise complex structures on $X$, if
$3r>4g-11$ where $g$ is the genus of the fibers of $f\co X\to
S^2$, then $i\co X_{2g-2-r}(f)\to P_r (f)$ is an embedding.
Further, $AJ\co X_r (f)\to P_r (f)$ restricts to
$AJ^{-1}(i(X_{2g-2-r}(f)))$ as a $\mathbb{P}^{r-g+1}$-bundle, and
is a $\mathbb{P}^{r-g}$-bundle over the complement of
$i(X_{2g-2-r}(f))$.
\end{lemma}
The reason for this is that in general
$AJ^{-1}(\mathcal{L})=\mathbb{P}H^0(\mathcal{L})$, which by
Riemann-Roch is a projective space of dimension
$r-g+h^1(\mathcal{L})$. The result of \cite{EH} ensures that for
$r>(4g-11)/3$ and for generic families of complex structures on
the $\Sigma_t$, none of the fibers of $f$ admit any line bundles
$\mathcal{L}$ with degree $r$ and $h^1(\mathcal{L})>1$; then
$Im(i)\subset P_r(f)$ consists of those bundles for which
$h^1(\mathcal{L})=h^0 (\kappa\otimes \mathcal{L}^{\vee})=1$.  To
see the bundle structure, rather than just set-theoretically
identifying the fibers, note that on any $\Sigma_t$, when we
identify the tangent space to $Pic^r\Sigma_t$ with
$H^0(\kappa_{\Sigma_t})$, the orthogonal complement of the
linearization $(AJ_*)_D$ at $D\in S^r\Sigma_t$ consists of those
elements of $H^0(\kappa_{\Sigma_t})$ which vanish along $D$ (this
follows immediately from the fact that, after choosing a basepoint
$p_0\in \Sigma_t$ and a basis $\{\phi_1,\ldots,\phi_g\}$ for
$H^0(\kappa_{\Sigma_t})$  in order to identify $Pic^{r}(\Sigma_t)$
with $\mathbb{C}^g/H^1(\Sigma_t,\mathbb{Z})$, $AJ$ is given by
$AJ(\sum
p_i)=\left(\sum\int_{p_0}^{p_i}\phi_1,\ldots,\sum\int_{p_0}^{p_i}\phi_g\right)$).
If $AJ(D)\notin Im (i)$, so that $H^0(\kappa-D)=0$, this shows
that $(AJ_*)_D$ is surjective, so that $AJ$ is indeed a submersion
away from $AJ^{-1}(Im\, i)$.  Meanwhile, if $\mathcal{L}=i(D')\in
Im (i)$, the above description shows that the only directions in
the orthogonal complement of any $Im(AJ_*)_D$ with
$AJ(D)=\mathcal{L}$ are those 1-forms which vanish at $D$, but
since $AJ(D)=i(D')$ such 1-forms also vanish at $D'$ and so are
also orthogonal to $Im(i_*)_{D'}$.  So if $AJ(D)=i(D')$, $Im
(AJ_*)_D$ contains $T_{i(D')}(Im\, i)$, implying that $AJ$ does in
fact restrict to $AJ^{-1}(Im\, i)$ as a submersion and hence as a
$\mathbb{P}^{r-g+1}$ bundle.

Smith's duality theorem, and also the vanishing result in this
paper, depend on the construction of almost complex structures
which are especially well-behaved with respect to the Abel-Jacobi
map.  From now on, we will fix complex structures on the fibers of
$X$ satisfying the conditions of Lemma \ref{BN}; these induce
complex structures on the fibers of the $X_r(f)$ and $P_r(f)$, but
on all of our spaces (including $X$) we still have the freedom to
vary the ``horizontal-to-vertical'' parts of the almost complex
structures.  Almost complex structures agreeing with these fixed
structures on the fibers will be called ``compatible.''

The following is established in the discussion leading to
Definition 6.4 of \cite{Smith}.

\begin{lemma} (\cite{Smith}) \label{dual} In the situation of Lemma \ref{BN}, for any compatible almost complex structure $J_1$ on $X_{2g-2-r}(f)$ and any
compatible $J_2$ on $P_{r}(f)$ such that $J_2|_{T(Im\,i)}=i_*J_1$,
there exist compatible almost complex structures $J$ on $X_r (f)$
with respect to which $AJ\co X_r(f)\to P_r(f)$ is
$(J,J_2)$-holomorphic. \end{lemma}

We outline the construction of $J$: Since $AJ\co AJ^{-1}(Im
\,i)\to X_{2g-2-r}(f)$ is a $\mathbb{P}^{r-g+1}$-bundle, given the
natural complex structure on $\mathbb{P}^{r-g+1}$ and the
structure $J_1$, the structures on $ AJ^{-1}(Im \,i)$  making this
fibration pseudoholomorphic correspond precisely to connections on
the bundle; since this bundle is the projectivization of the
vector bundle with fiber $H^0(\kappa-D)$ over $D$, a suitable
connection on the latter gives rise to a connection on our
projective-space bundle and thence to an almost complex structure
$J$ on  $AJ^{-1}(Im \,i)$ making the restriction of $AJ$
pseudoholomorphic.

To extend $J$ to all of $X_r (f)$, we first use the fact that, as
in Lemma 3.4 of \cite{DS}, \[ AJ_*\co \left(N_{ AJ^{-1}(Im
\,i)}X_r (f)\right)|_{AJ^{-1}(i(D))}\to (N_{Im\,i}P_r
(f))_{i(D)}\] is modeled by the map
\begin{align} \{(\theta,[x])\in V^*\times\mathbb{P}(V)|\theta(x)=0\}&\to V^* \nonumber\\
(\theta,[x])&\mapsto \theta, \nonumber \end{align}
   where $V=H^0(\kappa_{\Sigma_t}-D)$, so that the construction of Lemma 5.4 of \cite{DS} lets us extend $J$ to the closure of some open neighborhood
$U$ of $ AJ^{-1}(Im \,i)$.  But then since $AJ$ is a
$\mathbb{P}^{r-g}$-bundle over the complement of $ AJ^{-1}(Im
\,i)$, the problem of extending $J$ suitably to all of $X_r(f)$
amounts to the problem of extending the connection induced by $J$
from $\partial U$ to the entire bundle, which is possible because,
again, our bundle is the projectivization of a vector bundle and
connections on vector bundles can always be extended from closed
subsets.

Our vanishing results are consequences of the following:
\begin{lemma}(\cite{Smith},p.965)\label{neg}  Assume that $b^+(X)>b_1(X)+1$.  For any fixed compatible smooth almost complex structure $J_1$ on $X_{2g-2-r}(f)$ and for
generic smooth compatible almost complex structures $J_2$ such
that $J_2|_{Im\,i}=i_*J_1$, all $J_1$-holomorphic sections of
$P_r(f)$ are contained in $i(X_{2g-2-r}(f))$.
\end{lemma}

This follows from the fact that, as Smith has shown, the index of
the $\bar{\partial}$-operator on sections of $P_r(f)$ is
$1+b_1-b^+$, which under our assumption is negative, and so since
$J_2$ may be modified as we please away from $Im\,i$, standard
arguments show that for generic $J_2$ as in the statement of the
lemma all sections will be contained in $Im\,i$.

\section{Proof of Theorem \ref{famvan}} \label{vanproof}

\begin{lemma} \label{bigvan2} If $b^+(X)>b_1(X)+4n+1$, then
$\mathcal{FDS}^{n}_{f}(\alpha-2\sum e_i)=0$ for all $\alpha \in
H^2(X;\mathbb{Z})$ such that $r=\langle \alpha,[\Phi]\rangle$
satisfies $r>\max\{g(\Phi)+3n,2g(\Phi)-2\}$.\end{lemma}
\begin{proof}

Let $(J'_b)_{b\in X_n}$ be a smooth family of almost complex
structures on the relative Picard scheme $P_r(f)$ such that
\begin{itemize} \item[(i)] For each $b$, the map $G\co P_r(f)\to
S^2$ is pseudoholomorphic with respect to $J'_b $, and for all
critical values $t$ of $f$ $J$ agrees near $G^{-1}(t)$ with the
standard complex structure on the relative Picard scheme induced
by an integrable complex structure near $f^{-1}(t)$; \item[(ii)]
For each $b=(p_1,\ldots,p_n)$, where
$t_i=f\circ\pi^{(p_1,\ldots,p_{i-1})}(p_i)$, $J'_b$ also agrees
near each $G^{-1}(t_i)$ with the standard complex structure
induced by an integrable complex structure near $f^{-1}(t_i)$.
\end{itemize}  Thanks to the assumption that $b^+(X)>b_1(X)+4n+1$
and the fact that the index of the $\bar{\partial}$-operator on
sections of $P_r(f)$ is $1+b_1-b^+$, for a generic such family
$(J'_b)_{b\in X_n}$ there will be no $J'_b$ holomorphic sections
of $P_r(f)$ for any $b$.  Now, as in Section \ref{review}, since
$r>2g-2$, so that $AJ\co X_r(f)\to P_r(f)$ is a projective-space
bundle, we can construct a family $J_b$ of almost complex
structures on $X_r(f)$ such that $AJ\co X_r(f)\to P_r(f)$ is
$(J_b,J'_b)$-holomorphic for each $b$. By construction, for each
$b$ $J_b$ agrees with the standard complex structure on the
relative Hilbert scheme $F\co X_r(f)\to S^2$ near each singular
fiber and also near each $F^{-1}(t_i)$ where the $t_i$ are as
above.  Since $X^b$ is formed from $X$ by performing blowups at
points in $f^{-1}(t_i)$, for $b\in X'_n$ $J_b$ lifts to an almost
complex structure $\tilde{J}_b$ on $X^{b}_{r}(f^b)$ such that the
map $\Pi^b\co X^{b}_{r}(f^b)\to X_r(f)$ induced by blowup is
$(\tilde{J}_b,J_b)$-holomorphic.

Let $J^m$ be almost complex structures on the
$\mathcal{X}^{n}_{r}(f)$ from the Baire set in the definition of
$\mathcal{FDS}$ which converge to an almost complex structure that
agrees on each $X^{b}_{r}(f^b)$ with $\tilde{J}_b$.  If the
invariant were nonzero, we would obtain $J^m$-holomorphic sections
$s^m$ of some $X^{b_m}_{r}(f^{b_m})$ ($b_m\in X'_n$); after
passing to a subsequence we assume $b_m\to \bar{b}\in X_n$ (since
$X_n$, though not $X'_n$, is compact). By the definition of our
class of almost complex structures (see the text before Lemma
\ref{cpct2}) there are compatible almost complex structures
$J^{m}_{b_m}$ on $X_{r}(f)$ such that $\Pi^{b_m}\co
X^{b_m}_{r}(f^{b_m})\to X_r(f)$ is
$(J^m,J^{m}_{b_m})$-holomorphic; further, we will have
$J^{m}_{b_m}\to J_{\bar{b}}$.  So the $\Pi^{b_m}\circ s^m$ are
$J^{m}_{b_m}$-holomorphic sections of $X_{r}(f)$, whence after
passing to a subsequence they converge modulo bubbling to a
$J_{\bar{b}}$-holomorphic section $\bar{s}$. (As usual, even if
bubbling occurs, the bubble tree will contain a component which is
a $J_{\bar{b}}$-holomorphic section by virtue of the fact that all
bubbles will be contained in the fibers.)  But then $AJ\circ \bar{s}$
would be a $J'_{\bar{b}}$-holomorphic section, contradicting the
fact that no $J'_b$-holomorphic sections exist for any $b\in X_n$.
\end{proof}

The intermediate case where
$\max\{g(\Phi)+3n+d(\alpha),(4g(\Phi)-11)/3\}<r\leq 2g(\Phi)-2$
takes slightly more work.  In this case, as in Section
\ref{review} we use the fact  that combining the Abel-Jacobi map
with Serre duality gives a map
\[ i\co X_{2g-2-r}(f)\to P_{r}(f);\]
as before since $3r>4g-11$ generic choices of the complex
structures on the fibers of $f$ result in this map being an
embedding.  Similarly to the proof of Lemma \ref{bigvan2},
consider families of almost complex structures $J''_b$ ($b\in
X_n$) on $X_{2g-2-r}(f)$ which make $X_{2g-2-r}(f)\to S^2$
holomorphic and are standard near the singular fibers and near the
fibers containing the points which are blown up to form $X^b$.
Form almost complex structures $J'_b$ on $P_r(f)$ restricting to
$i(X_{2g-2-r}(f))$ as $i_*J''_b$ and which are also standard near the
singular fibers and near the fibers containing the points which
are blown up to form $X^b$.  The fact that $b^+>b_1+1+4n$ implies
that if the family $J'_b$ is chosen generically among almost
complex structures with this property, then any $J'_b$ holomorphic
sections of $P_r(f)$ for any $b$ must be contained in
$i(X_{2g-2-r}(f))$.

We then form almost complex structures $J_b$ on $X_r(f)$ such that
$AJ\co X_r(f)\to P_r(f)$ is $(J_b,J'_b)$-holomorphic.  As in the
proof of Lemma \ref{bigvan2}, a nonvanishing invariant
$\mathcal{FDS}^{n}_{f}(\alpha-2\sum e_i)$ would give rise to a
sequence of sections of $X_r(f)$ in the homotopy class
$c_{\alpha}$ which converge modulo bubbling to a
$J_{\bar{b}}$-holomorphic section $\bar{s}$ of $X_r(f)$.  Since
all fibers of $f$ are irreducible, any bubbles that arise will
descend to a multiple covering of one of the fibers of $f$, and so
for some $m\geq 0$ we will have $\bar{s}\in c_{\alpha-m PD[\Phi]}$
where as usual $[\Phi]$ is the class of the fiber.

$AJ\circ \bar{s}$ will then be a $J'_{\bar{b}}$-holomorphic
section of $P_r(f)$, and so must be contained in
$i(X_{2g-2-r}(f))$.  By the construction of $i$, then,
$i^{-1}\circ AJ\circ \bar{s}$ is a $J''_{\bar{b}}$-holomorphic
section of $X_{2g-2-r}(f)$ in the homotopy class
$c_{\kappa_X-\alpha+mPD[\Phi]}$.

Now one computes using the adjunction formula for the fiber $\Phi$
that
\begin{align*}
d(\kappa_{X^b}-\alpha+mPD[\Phi])
&=d(\kappa_X-\alpha)+d(m\Phi)+m\langle\kappa_X-\alpha,[\phi]\rangle
\\ &=d(\alpha)-\frac{m}{2}\langle \kappa_X,[\Phi]\rangle+m\langle
\kappa_X-\alpha,[\Phi]\rangle \\ &=d(\alpha)-m(r-g(\Phi)+1).
\end{align*}

Thus by choosing the $4n$-real-dimensional family $J''_b$
generically we ensure that $m=0$ thanks to the assumption that
$r>d(\alpha)+g(\Phi)+3n$ in the statement of the theorem.

Now take a family of almost complex structures $j_b$ on $X$ which
are standard near the singular fibers of the fibrations $f$ and
also near the fibers containing the points blown up to form $X^b$;
these induce tautological almost complex structures
$\mathbb{J}_{j_b}$ on $X_{2g-2-r}(f)$. Let $J''^{m}_{b}$ be
families of smooth almost complex structures on $X_{2g-2-r}(f)$
which are generic in the sense of the previous paragraph and which
converge in H\"older norm to the $\mathbb{J}_{j_b}$. For each $m$
there is some $b_m$ such that $J''^{m}_{b_m}$ admits a holomorphic
section in the class $c_{\kappa_{X}-\alpha}$, so Gromov
compactness guarantees the existence of a
$\mathbb{J}_{j_{b_0}}$-holomorphic section of some $X_{2g-2-r}(f)$
in $c_{\kappa_{X}-\alpha}$ for some $b_0$; this section then
tautologcally corresponds to a $j_{b_0}$-holomorphic curve $C$
Poincar\'e dual to $\kappa_{X}-\alpha$; setting $j=j_{b_0}$, this
is the curve that we desire.

To get the $j$-holomorphic curve Poincar\'e dual to $\alpha$, we
simply consider the almost complex structures $j_b$ on the members
$X^b$ of the family blowup induced in the almost complex category
by $j$.  Let $j_{b}^{m}$ be a sequence of almost complex
structures $C^0$-approximating the $j_b$ which are integrable near
the exceptional spheres, and apply Gromov compactness to a
sequences of almost complex structures on $\mathcal{X}^{n}_{r}(f)$
whose restrictions to $X^{b}_{r}(f^b)$ H\"older-approximate the
family $\mathbb{J}_{j_{b}^{m}}$; in this way our nonvanishing
invariant guarantees the existence of a
$\mathbb{J}_{j_{b_m}^{m}}$-holomorphic section of some
$X^{b_m}_{r}(f^{b_m})$ in the class $c_{\alpha-2\sum e_i}$ and so
of a $j_{b_m}^{m}$-holomorphic curve Poincar\'e dual to
$\alpha-2\sum e_i$.  Appealing to Gromov compactness for these
curves then gives a $j_b$-holomorphic curve, and this latter is
sent by the blowdown map to the $j$-holomorphic curve which we
desire. Theorem \ref{famvan} is thus proven.

If $X$ admits an integrable complex structure $j$ making the
fibration holomorphic, then for our original family of almost
complex structures $j_b$ we can take the constant family $j$,
justifying a statement made near the end of the introduction.  For
arbitrary $j$, though, this argument does not work, because it was
crucial in the construction of the curve Poincar\'e dual to
$\kappa_{X}-\alpha$ that each of the $j_b$ was integrable near the
fibers containing the points blown up in forming $X^b$.

\section{Two technical matters}

\subsection{Blowing up a point in an almost complex
manifold} \label{app2}In the proof of Theorem \ref{famsame} we
have used the fact that, if $\pi\co X'\to X$ is the blowup of a
4-manifold at a point and $J$ is an almost complex structure on
$X$, then there is a Lipschitz almost complex structure $J'$ on
$X'$ such that $\pi$ is $(J',J)$-holomorphic.  Since we have not
found a proof of this fact in the literature, we present one here.
As the dimension of $X$ does not affect the argument, we prove the
result for almost complex manifolds of arbitrary complex dimension
$n$. The blowup, of course, has the effect of replacing the point
$p$ being blown up with an exceptional divisor $E\cong
\mathbb{C}P^{n-1}$; we note that, as will be seen in the proof,
$J'|_{TE}$ agrees with the standard complex structure on
$\mathbb{C}P^{n-1}$.  If $(X,\omega)$ is symplectic, recall from,
\emph{e.g.}, Chapter 7 of \cite{MS2} that $X'$ can be endowed with
symplectic forms $\omega_{\ep}$ for small $\ep>0$, with the
parameter $\ep$ reflecting the size of the exceptional divisor $E$
in the symplectic manifold $(X',\omega_{\ep})$.  One can easily
check that if the almost complex structure $J$ on $X$ is
$\omega$-tame, then $J'$ will be $\omega_{\ep}$-tame for small
enough $\ep$.

 Our method only
proves Lipschitz regularity for $J'$; it is unclear whether $J'$
is differentiable in directions normal to $E$.  In principle, one
would also like to be able to blow up almost complex submanifolds
$V\subset (X,J)$ of arbitrary dimension in the almost complex
category.  Our method does not readily extend to show that the
pullback of $J$ under the blowup extends even continuously over
the exceptional divisor of the blowup when $\dim V>0$.
Nonetheless, the case of blowing up a point suffices for our
application.

We begin with the following lemma, which will later be used to
construct coordinate charts on the blowup.

\begin{lemma} \label{fix} Let $J$ be an almost complex structure on
$\mathbb{C}^n$ agreeing at the origin with the standard complex
structure $J_0$.  Given $\kappa_0\in \mathbb{C}P^{n-1}$ there
exists a constant $\rho_0$ with the following property. Let
$\rho<\rho_0$ and let $U_{\rho}$ be the ball of radius $\rho$
around $\kappa_0$ in $\mathbb{C}P^{n-1}$ and $D_{\rho}$ the disc
of radius $\rho$ in $\mathbb{C}$.  There is a smooth map \[
\Theta\co D_{\rho}\times U_{\rho}\to \mathbb{C}^n \] such that
each $\Theta|_{D_{\rho}\times \{\kappa\}}$ ($\kappa\in
U_{\rho}\subset \mathbb{C}P^{n-1}$) is an embedding whose image is
a $J$-holomorphic disc which is tangent at the origin to the line
$l_{\kappa}\subset \mathbb{C}^n$ determined by $\kappa$.
  \end{lemma}
\begin{proof}  The proof quite closely parallels some of the
arguments in Section 5 of \cite{TSW}; we outline it for
completeness. By a complex linear change of coordinates we may
assume that $\kappa_0 =[1:0:\cdots :0]$. Where
$c=(c_1,\ldots,c_{n-1})\in (D_{\rho})^{n-1}$ and
$\kappa=[1:\kappa_1:\cdots :\kappa_{n-1}]$ is close to
$[1:0:\cdots :0]$, we search for a $J$-holomorphic disc
\[ q_{c,\kappa}(z)=(z,c_1+\kappa_1
z+u_1(c,\kappa,z),\ldots,c_{n-1}+\kappa_{n-1}z+u_{n-1}(c,\kappa,z))
\] defined for $z\in D_{\rho}$.  As in \cite{TSW}, this is equivalent to a system of equations \[ \frac{\partial
u_i}{\partial
\bar{z}}=Q_i\left(c,\kappa,u_1(c,\kappa,z),\ldots,u_{n-1}(c,\kappa,z)\right)
\] such that for certain constants $\gamma_k$ we have
\begin{equation}\label{qsmall} \|Q_i\|_{C^k}\leq \gamma_k\|J-J_0\|_{C^k(D_{2\rho}^{n})}.\end{equation}  Note that by
decreasing $\rho$ and rescaling the coordinates we can make the
right hand side of (\ref{qsmall}) as small as we like.

Now introduce a cutoff function $\chi_{\rho}\co\mathbb{C}\to
[0,1]$ which equals $1$ for $|z|<\rho$ and $0$ for $|z|>3\rho/2$,
and search for a solution to \[ \frac{\partial
u_i}{\partial\bar{z}}=\chi_{\rho}Q_i \quad (i=1,\ldots,n-1)\] by,
on the class of $(n-1)$-tuples of $C^{2,1/2}$ functions $u_i$
restricting to the circle of radius $4\rho$ around zero in the
span of $\{e^{ik\theta}|k<0\}$, searching for a tuple
$(u_1,\ldots,u_{n-1})$ obeying \begin{equation}\label{fixedpt}
\left(
u_i(z)=\frac{1}{\pi}\int\frac{\chi_{\rho}Q_i(c,\kappa,u_i(c,z))}{z-w}d^2w\right)_{i=1,\ldots,n-1}
\end{equation}  Applying the contractive mapping theorem on this
class of functions (viewed as a Banach space using the
$(n-1)$-fold direct sum of the norm used on p. 886 of \cite{TSW}),
thanks to the smallness of the $Q_i$ we can find a unique small
solution of (\ref{fixedpt}).  Furthermore as in Lemma 5.5 of
\cite{TSW} the solution varies smoothly in each of $z$, $c$, and
$\kappa$,  and satisfies bounds \[ \left|\frac{\partial
u}{\partial c_i}\right|<C\rho,\quad \left|\frac{\partial
u}{\partial \kappa_i}\right|<C\rho^2,\quad
\|u\|_{C^0}<C(\rho^2+\rho(|c|+|\kappa|)),\quad
\|u\|_{C^1}<C(\rho+(|c|+|\kappa|)).\]

Letting $\sigma$ denote the map which assigns to $(c,\kappa)$ the
pair consisting of $q_{c,\kappa}(0)$ and the tangent space to
$Im\, q_{c,\kappa}$ at $q_{c,\kappa}(0)$, the implicit function
theorem then allows us to solve the equation
$\sigma(c,\tilde{\kappa})=((0,\ldots,0),\kappa)$ for $c$ and
$\tilde{\kappa}$ in terms of $\kappa$. The desired map $\Theta$ is
then \begin{align*} \Theta\co D_{\rho}\times U_{\rho}&\to \mathbb{C}^{n} \\
(z,\kappa)&\mapsto q_{c(\kappa),\tilde{\kappa}(\kappa)}(z).
\end{align*}
\end{proof}

For any even-dimensional manifold $X$ with $p\in X$, we  form the
blowup $X'$ of $X$ at $p$ as a topological manifold by removing a
ball $B^{2n}$ around $p$, embedding $B^{2n}$ in $\mathbb{C}^n$ in
standard fashion, and replacing $B^{2n}$ in $X$ by $B'= \{(l,e)\in
\mathbb{C}P^{n-1}\times \mathbb{C}^n|e\in l\cap B^{2n}\}$. The
blowdown map $\pi\co X'\to X$ is of course just the identity
outside $B'$ and the map $(l,e)\to e$ inside $B'$.  The
exceptional divisor is $E=\{(l,e)\in B'|e=0\}\subset X'$.

If $\kappa_0=[1:0:\cdots :0]$ in Lemma \ref{fix} and we write
$\kappa$ near $\kappa_0$ as $[1:\kappa_1\cdots:\kappa_{n-1}]$, the
map $\Theta$ has the form
\[ (z,\kappa)\mapsto (z,\kappa_1z+
\tilde{u}_1(\kappa,z),\ldots,\kappa_{n-1}z+\tilde{u}_{n-1}(\kappa,z))
\] where the $\tilde{u}_i$ are smooth functions satisfying
$|\tilde{u}_i(\kappa,z)|<C|z|^2$ for an appropriate constant $C$.
(In the notation of the proof of Lemma \ref{fix},
$\tilde{u}_i(\kappa,z)=u_i(c(\kappa),\tilde{\kappa}(\kappa),z)+(\tilde{\kappa}_i(\kappa)-\kappa_i)z$.)

We hence obtain a local homeomorphism
$\tilde{\Theta}=\tilde{\Theta}_{\kappa_0}\co D^2\times D^2\to
\tilde{\mathbb{C}}^n$ such that, where $\pi\co
\tilde{\mathbb{C}}^n\to\mathbb{C}^n$ is the blowdown, $\pi\circ
\tilde{\Theta}=\Theta$.  We use the $\tilde{\Theta}_{\kappa_0}$ as
$\kappa_0$ varies over $\mathbb{C}P^{n-1}$ as an atlas for
$\tilde{\mathbb{C}}^n$ near the exceptional divisor $E$ (away from
$E$ we of course just use charts pulled back by $\pi$ from charts
on $\mathbb{C}^n$ not containing the origin).  From the definition
of the $\tilde{\Theta}_{\kappa_0}$ and the fact that tangencies of
$J$-holomorphic curves in $\mathbb{C}^n$ are $C^1$-diffeomorphic
to tangencies between $J_0$-holomorphic curves \cite{Sik}, one can
see that the transition functions have the form \[
\tilde{\Theta}_{\kappa_0}^{-1}\circ\tilde{\Theta}_{\kappa'_0}(z,\kappa_1,\ldots,\kappa_{n-1})=\qquad
\qquad \]
\[\left(z,\kappa_1+z^{-1}(f_1(\kappa)z^2+O(|z|^3)),\ldots,\kappa_{n-1}+z^{-1}(f_{n-1}(\kappa)z^2+O(|z|^3))\right),\]
and in particular are $C^1$.  We have thus provided an atlas for
$\tilde{\mathbb{C}}^{n}$ as a $C^1$ manifold.

This atlas depends on the almost complex structure $J$, and it is
worth noting that the charts corresponding to different $J$ might
not be $C^1$-related.  For example, for a particular $J$
$\Theta_{[1:0:\cdots:0]}$ could conceivably have the form
\[\Theta_{[1:0:\cdots:0]}(z,\kappa_1,\ldots,\kappa_n)=(z,\kappa_1z+\bar{z}^2,\kappa_2z,\ldots,\kappa_{n-1}z).\]
In this case, in terms of the \emph{standard} smooth coordinates
on $\tilde{\mathbb{C}}^n$ (equivalently, those induced by the
above construction using the standard complex structure $J_0$ ),
\[\tilde{\Theta}_{[1:0:\cdots:0]}(z,\kappa_1,\ldots,\kappa_{n-1})=(z,\kappa_1+\bar{z}^2/z,\kappa_2,\ldots,\kappa_{n-1}),\]
which  is Lipschitz but not $C^1$ along the exceptional divisor
$\{z=0\}$.  Of course, these resulting manifolds are still
abstractly $C^1$-diffeomorphic; this is somewhat reminiscent of
the fact that distinct complex structures on a Riemann surface
$\Sigma$ induce smooth charts on the symmetric products
$S^d\Sigma$ which are related by transition maps that are only
Lipschitz, as noted for instance in Remark 4.4 of \cite{Sa}.

\begin{prop}  Let $\pi\co \tilde{\mathbb{C}}^n\to \mathbb{C}^n$
denote the blowup of $\mathbb{C}^n$ at the origin, and let $J$ be
an almost complex structure on $\mathbb{C}^n$ agreeing with the
standard almost complex structure $J_0$ at the origin.  Then there
is a unique Lipschitz continuous almost complex structure
$\tilde{J}$ on $\tilde{\mathbb{C}}^n$ such that $\pi$ is a
$(\tilde{J},J)$ holomorphic map. \end{prop}
\begin{proof} Let $E\cong\mathbb{C}P^{n-1}$ denote the exceptional
divisor of the blowup $\pi$.  Of course, $\pi$ restricts to a
diffeomorphism $\tilde{\mathbb{C}^n}\setminus E\to
\mathbb{C}^n\setminus (0,\ldots,0)$, so our $\tilde{J}$ must agree
away from $E$ with $\pi^{*}J=\pi_{*}^{-1}\circ J\circ\pi_*$ away
from $E$ and uniqueness even of a continuous almost complex
structure $\tilde{J}$ is clear from the fact that
$\tilde{\mathbb{C}}^{n}\setminus E$ is dense in
$\tilde{\mathbb{C}}^{n}$.  We show now that $\pi^{*}J$ extends
over $E$  in Lipschitz fashion by exhibiting a Lipschitz
continuous basis of vector fields for its antiholomorphic tangent
space $T^{0,1}\subset T\tilde{\mathbb{C}}^n\otimes \mathbb{C}$
near any given point $x\in E$.

Lemma \ref{fix} and the remarks thereafter provide us with one
element of this basis: the maps $\Theta_{\kappa_0}$ map each
$D_{\rho}\times \{\kappa\}$ diffeomorphically to a $J$-holomorphic
disc $\Delta_{\kappa}$ in $\mathbb{C}^n$ in a way that varies
smoothly in $\kappa$.  We then obtain a (complexified) vector
field $\tilde{\alpha}_{\kappa}$ along each $D_{\rho}\times
\{\kappa\}$ defined by the property that
$\alpha_{\kappa}=(\Theta_{\kappa_0})_{*}\tilde{\alpha}$ generates
the $J$-antiholomorphic tangent space to $\Delta_{\kappa}$.
Choosing the $\alpha_{\kappa}$ to depend smoothly on $\kappa$
causes the $\tilde{\alpha}_{\kappa}$ to do so as well, and so to
give a vector field $\alpha$ on a neighborhood of our basepoint
$x$ which is transverse to $E$ and which is antiholomorphic for
the pulled back almost complex structure $\pi^{*}J$ where the
latter is defined.

After a complex linear change of coordinates on $\mathbb{C}^n$ we
may assume that $x=([1:0:\cdots:0],(0,\ldots,0))$ and
$\pi(x)=(0,\ldots,0)$. In terms of the coordinate chart given by
$\Theta_{[1:0\cdots:0]}$, the blowdown map $\pi$ has the form
\[ (s,t_1,\ldots,t_{n-1})\mapsto (s,st_1+u_1(s,t_1,\ldots,t_{n-1}),\ldots,st_{n-1}+u_{n-1}(s,t_1,\ldots,t_{n-1})),\] where
$|u_{i}(s,t_1,\ldots,t_{n-1})|<C|s|^2$.  Away from the exceptional
sphere $s=0$, this is a diffeomorphism whose complexified
linearization with respect to the coordinates
$(s,\bar{s},t_1,\bar{t}_1,\ldots,t_{n-1}\bar{t}_{n-1})$ has
inverse of the form
\[
((\pi_*)^{-1})_{\pi(s,t_1,\ldots,t_{n-1})}=\left(\begin{matrix}
1&0&\cdot&\cdot&\cdot&\cdot&0\\
0&1&0&\cdot&\cdot&\cdot&0\\
-t_1/s&0&1/s&0&\cdot&\cdot&0\\
0&-\bar{t}_1/\bar{s}&0&1/\bar{s}&0&\cdots&0\\
\vdots&\vdots& \vdots& \ddots&\ddots&\ddots&\vdots\\
-t_{n-1}/s&0&\cdot&\cdots&0&1/s&0\\
0&-\bar{t}_{n-1}/s&0&\cdot&\cdots&0&
1/\bar{s}\end{matrix}\right)+B(s,t_1,\ldots,t_n),\] where $B$ is
smooth away from $s=0$ and bounded (but not necessarily
continuous) as $s\to 0$.

 Write the coordinates on $\mathbb{C}^n$ as
$(w,z_1,\ldots,z_{n-1})$. Since $J$ agrees with $J_0$ at the
origin, for $i=1,\ldots, n-1$ there are $J$-antiholomorphic vector
fields \begin{align*} \beta_i=\partial_{\bar{z}_i}+\sum_j
a_{ij}(z_1,\ldots,z_n)\partial_{z_j}+\sum_{j\neq
i}b_{ij}(z_0,\ldots,z_n)\partial_{\bar{z}_j}+c_i(z_0,\ldots,z_n)\partial_{w},\end{align*}
where $a_{ij}(0,\ldots,0)=b_{ij}(0,\ldots,0)=c_i(0,\ldots,0)=0$.
Away from $E$ we then have \begin{align*}
(\pi_{*}^{-1}\beta_i)_{(u,v_1,\ldots,v_{n-1})}&=\frac{1}{\bar{u}}\partial_{\bar{v}_i}+\sum_j
a_{ij}(\pi(u,v_1,\ldots,v_{n-1}))\left(\frac{1}{u}\partial_{v_j}\right)\\&+
\sum_{j\neq i}
b_{ij}(\pi(u,v_1,\ldots,v_{n-1}))\left(\frac{1}{\bar{u}}\partial_{\bar{v}_j}\right)\\&+c_i(\pi(u,v_1,\ldots,v_{n-1}))\left(\partial_u-\sum_j\frac{v_j}{u}\partial_{v_j}\right)+\tilde{\gamma_i}\end{align*}
where $\tilde{\gamma_i}=B\beta_i$ has bounded coeffecients.  So
\begin{align*}
\tilde{\beta}_i&:=\bar{u}\pi_{*}^{-1}\beta_i=\partial_{\bar{v}_i}+\sum\frac{\bar{u}}{u}\left(a_{ij}(u,uv_1,\ldots,uv_{n-1})-v_jc_{ij}(u,uv_1,\ldots,uv_{n-1})\right)\partial_{v_j}
\\&+\sum_{j\neq
i}b_{ij}(u,uv_1,\ldots,uv_{n-1})\partial_{\bar{v}_j}+\bar{u}c_i(u,uv_1,\ldots,uv_{n-1})\partial_w+\bar{u}\tilde{\gamma_i}\end{align*}
is an antiholomorphic tangent vector for $\pi^{*}J$ away from
$E=\{u=0\}$. Further, we note that since $a_{ij},b_{ij}$, and
$c_i$ are differentiable and vanish at the origin while
$\tilde{\gamma_i}$ is $L^{\infty}$, so that
$|a_{ij}(u,uv_1,\ldots,uv_{n-1})|$,
$|b_{ij}(u,uv_1,\ldots,uv_{n-1})|$,
$|c_{i}(u,uv_1,\ldots,uv_{n-1})|$, and $\|\bar{u}\gamma_i\|$ are
all bounded by a constant times $|u|$, $\tilde{\beta}_i$ extends
over $E$ in Lipschitz fashion, agreeing with
$\partial_{\bar{v}_i}$ at $E$.

Hence, defining $\tilde{J}$ near $x$ by \[
T^{0,1}_{\tilde{J}}=\langle
\tilde{\alpha},\tilde{\beta_1},\ldots,\tilde{\beta}_{n-1}\rangle,\]
we see that $\tilde{J}$ is Lipschitz and agrees with $\pi^{*}J$
where the latter is defined.  So since $\tilde{J}$ preserves $TE$
and since at each point of $E$ there is a $\tilde{J}$-holomorphic
disc transverse to $E$ mapped holomorphically to a $J$-holomorphic
disc by $\pi$, we conclude that $\pi\co \tilde{\mathbb{C}}^n\to
\mathbb{C}^n$ is $(J',J)$-holomorphic.\end{proof}

\begin{cor}  Let $(X,J)$ be an almost complex manifold with $p\in
X$, and let $X'$ denote the blowup of $X$ at $p$.  Then there is a
unique almost complex structure $J'$ on $X'$ which is Lipschitz
continuous such that $\pi\co X'\to X$ is $(J',J)$-holomorphic.
Further $J'$ restricts to $E$ as the standard complex structure on
$\mathbb{C}P^{n-1}$. \end{cor}
\begin{proof}  Since $\pi$ is a diffeomorphism away from $E=\pi^{-1}(p)$ (which thus determines $J'$ on $X'\setminus E$ as the
smooth almost complex structure $\pi^{*}J$), this follows from the
proposition and its proof by choosing a chart around $p$ which
sends $(p,J|_{T_p X})$ to $(0,J_{0}|_{T_0\mathbb{C}^n})$ in
$\mathbb{C}^n$ (as may easily be done by modifying any chart
around $p$ by an appropriate real linear map).
\end{proof}

\subsection{The diagonal in the relative Hilbert
scheme} \label{app1}

Let $F\co \mathcal{H}_r\to D^2$ denote the $r$-fold relative
Hilbert scheme of the map $f\co (z,w)\mapsto zw$;  the spaces
$\mathcal{H}_r\times \mathbb{C}^{s-r}$ form the local models for
the relative Hilbert scheme $X_s(g)$ of a Lefschetz fibration $g$
near points of $X_s(g)$ which correspond to divisors containing
$r$ copies of a critical point of $g$.  In this subsection we
prove the fact, used in the proof of the compactness result
underlying the construction of $\mathcal{FDS}$, that at a point in
the diagonal $\Delta$ of the relative Hilbert scheme
$\mathcal{H}_r$ corresponding to  the divisor in the nodal fiber
$f^{-1}(0)$ consisting of $r$ copies of $(0,0)$, the tangent cone
to the diagonal is contained in the tangent cone to the fiber
$F^{-1}(0)\subset \mathcal{H}_r$. (Note that since the natural map
$F^{-1}(t)\to S^r f^{-1}(t)$ is an isomorphism if and only if $t$
is a regular value of $f$, there are many points in $F^{-1}(0)$
corresponding to $\{(0,0),\ldots,(0,0)\}$, as will be seen later
on when we review the definition of $\mathcal{H}_r$.) Our proof of
this fact uses the description of the relative Hilbert scheme in
terms of linear algebra provided in Section 3 of \cite{Smith}
based on work of Nakajima \cite{N}, and boils down to a rather
arcane fact about the discriminants of the characteristic
polynomials of certain matrices.  It would certainly not surprise
us if there exists a more elegant way of proving this result via
algebraic geometry, but the argument we give presently is the only
one we have at the moment.  As will be seen later on, the relevant
characteristic polynomials have the form considered in the
following lemma.

\begin{lemma} \label{discrim}There is a universal, nonzero polynomial
$P(c_{k+1},\ldots ,c_{k+l+1})$ with $P(0,\ldots,0)=0$ such that,
given a degree $r=k+l+1$ polynomial \begin{equation}\label{poly}
f(x)=x^{r}+\sum_{a=1}^{k}\ep(c_a+O(\ep))x^{r-a}+\sum_{b=1}^{l+1}\ep^b(c_{k+b}+O(\ep))x^{l+1-b},
\end{equation}
the discriminant $\delta (f)$ of $f$ has the form
\begin{equation}\label{discform} \delta (f)=P(c_{k+1},\ldots
,c_{k+l+1})\ep^{r+l^2-1}+O(\ep^{r+l^2}).\end{equation} \end{lemma}

\begin{proof}
For $i=0,\ldots,r=k+l+1$, let $a_i$ be the coefficient of
$x^{r-i}$ in $f$ (so in particular $a_0=1$).  Recall that $\delta
(f)=(-1)^{r(r-1)/2}a_{0}^{-1} Res(f,f') $ (``$Res$'' denoting the
resultant; see, \emph{e.g.}, Section V.10 of \cite{Lang}), so it
suffices to prove the expansion (\ref{discform}) for $Res(f,f')$.
$Res(f,f')$ is given as the determinant \begin{equation}
\label{res} \left|
\begin{matrix}
a_0 & a_1 & a_2 &\cdots &\cdots & a_r   &     & &  \\
    & a_0   & a_1 & a_2 &\cdots & \cdots & a_r & &  \\
    &       & \cdot & \cdot & \cdot & \cdot & \cdot & \cdot & \\
    &       &       & a_0 & a_1 & a_2 &\cdots & \cdot & a_r \\
ra_0 & (r-1)a_1 & (r-2)a_2 &\cdots & a_{r-1}   &     & &  \\
    &\cdot &\cdot & \cdot  & \cdot & \cdot &   & \\
    &       & \cdot & \cdot   & \cdot  & \cdot  & \cdot & \\
    &       &       & ra_0 & (r-1)a_1 & (r-2)a_2 &\cdots & a_{r-1} & \\
    &       &       &     & ra_0 & (r-1)a_1 & (r-2)a_2 &\cdots & a_{r-1}

    \end{matrix}\right| \end{equation}

Each term in the expansion of this determinant will be a constant
times $\prod_{j=0}^{r} a_{j}^{i_j}$ for some natural numbers $i_j$
satisfying \[ \sum i_j = 2r-1 \] (since this is a $(2r-1)\times
(2r-1)$ matrix) and \[ \sum ji_j =r(r-1) \] (since if the roots of
$f$ are $\alpha_1,\ldots,\alpha_r$, the discriminant
$\prod_{a<b}(\alpha_a-\alpha_b)^2$ has degree $r(r-1)$ in the
$\alpha_b$, while the coefficient $a_j$ has degree $j$ in the
$\alpha_b$).  Let \[ e(i_0,\ldots,i_r)=\max\{e\in
\mathbb{N}|a_{0}^{i_0}\cdots a_{r}^{i_r}=O(\ep^e)\}.\]  To prove
the lemma we need to show that: \begin{itemize} \item[(i)]  For
each $\prod a_{j}^{i_j}$ appearing in the expansion of the
resultant (\ref{res}), $e(i_0,\ldots,i_r)\geq r+l^2-1$, with
equality implying that $i_1=\cdots =i_k=0$ (the latter condition
being needed to show that our polynomial $P$ depends only on
$c_{k+1},\ldots, c_r$ and vanishes when all of these $c_j$ are
$0$); and \item[(ii)] There are particular values of the $c_j$ for
which $Res(f,f')\neq O(\ep^{r+l^2})$. \end{itemize}

Point (ii) above is easy: in the statement of the lemma, let \[
c_j=\begin{cases} 1 & i=k+1,n \\ 0 &\text{otherwise,}
\end{cases}\] so that \[ f(x)=x^r+(\ep
+O(\ep^2))x^l+(\ep^{l+1}+O(\ep^{l+2})).\]  We then see that the
unique lowest-order term in the expansion  of the determinant
\ref{res} is obtained by choosing $a_0=1$ from the first $k+1$
columns, $(r-k-1)a_{k+1}=l\ep+O(\ep^2)$ from the next $r$ columns,
and $a_n=\ep^{l+1}+O(\ep^{l+2})$ from the last $l-1$ columns, so
that
\[ Res(f,f')=\pm(l\ep)^r(\ep^{l+1})^{l-1}+\mbox{higher order
terms}=\pm l^n\ep^{r+l^2-1}+O(\ep^{r+l^2}).\]

We now set about the proof of point (i).  Assume that $\prod
a_{j}^{i_j}$ is a term appearing in the expansion of the
determinant (\ref{res}).  Let $q$ be the quotient and $p$ be the
remainder when $\sum_{m=0}^{l}mi_{r-m}$ is divided by $l$, and set
$s=\sum_{m=0}^{l}i_{r-m}-q$ (note that the above sums only go up
to $l=r-k-1$).  We then have \[
\sum_{j=r-l}^{r}ji_j=\sum_{j=0}^{l}(r-l+j)i_{r-l+j}=(r-l)q+rs-p.\]

Now since $\sum_{j=0}^{r}i_j=\sum_{j=0}^{k}i_j+q+s=2r-1$ and since
$2r-1=r+k+l$, we see \begin{align}
s&=2r-1-q-\sum_{j=0}^{k}(i_j-1)-(k+1) \nonumber \\
&=r+l-1-q-\sum_{j=0}^{k}(i_j-1). \nonumber \end{align}

Hence \begin{align} r^2-r&=
\sum_{j=0}^{r}ji_j=\sum_{j=0}^{k}ji_j+q(r-l)+rs-p \nonumber \\
&=\sum_{j=0}^{k}ji_j+q(r-l)-p+r(r+l-1-q-\sum_{j=0}^{k}(i_j-1))
\nonumber \\ &=r^2-r+l(r-q)-p+\sum_{j=0}^{k}(ji_j-r(i_j-1)),
\nonumber
\end{align} \emph{i.e.},
\begin{equation}\label{lrq}
l(r-q)=p+\sum_{j=0}^{k}(r(i_j-1)-ji_j).
\end{equation}

Meanwhile \begin{align} e(i_0,\ldots,i_r)&=\sum_{j=1}^{k}i_j
+\sum_{j=0}^{l}(1+j)i_{r-l+j} \nonumber \\ &=\sum_{j=1}^{k}
i_j+q+s(l+1)-p \nonumber\\
&=\sum_{j=1}^{k}i_j+q+\left(r+l-1-q-\sum_{j=0}^{k}(i_j-1)\right)(l+1)-p
\nonumber
\\&=l(r-q)+r+l^2-1-(l+1)\sum_{j=0}^{k}(i_j-1)+\sum_{j=1}^{k}i_j-p
\nonumber \\
&=r+l^2-1+\sum_{j=0}^{k}\left(r(i_j-1)-ji_j\right)-(l+1)\sum_{j=0}^{k}(i_j-1)+\sum_{j=1}^{k}i_j
\nonumber \\ &=
r+l^2-1+k\sum_{j=0}^{k}(i_j-1)+\sum_{j=1}^{k}(1-j)i_j
\label{rearrange}, \end{align} where in the penultimate equality
we have used (\ref{lrq}) and in the last we have used the fact
that $r-(l+1)=k$.

In our term $\prod a_j^{i_j}$ in the expansion of the determinant
(\ref{res}), each of those $a_j$ which are chosen from the first
$(k+1)$ columns necessarily has $j\leq k$.  For each $j$ write
$i_j=w_j+z_j=w_j+x_j+y_j$ where $w_j$ denotes the number of
$a_j$'s chosen from the first $(k+1)$ columns and $x_j$ denotes
the number of $a_j$'s chosen from columns $k+2$ through $2k+1$;
evidently $w_j=0$ for $j>k$ while $\sum_{j=0}^{k} w_j=k+1$,
\emph{i.e.},
\begin{equation} \label{k+1} \sum_{j=0}^{k}(w_j-1)=0.
\end{equation}

Rearrange our term $\prod_{j=0}a_{j}^{i_j}$ as \[ a_{p_1}\cdots
a_{p_{2r-1}},\] where the entry $a_{p_n}$ is culled from the $nth$
column in the matrix in (\ref{res}); label the row from which
$a_{p_n}$ is taken as $m_n$. Denoting \[ \bar{m}=\begin{cases} m & m\leq r-1 \\
m+1-r & m\geq r \end{cases}, \] we see from the form of the
resultant matrix that
\[\bar{m}_n= n-p_n.\]  Consider the quantity \[ \sum_{n=1}^{2k+1}
\bar{m}_n.\]  Obviously, the way to minimize this quantity is by
using rows $1,2,\ldots,k,r,r+1,\ldots,r+k$ (or, just as well, rows
$1,\ldots,k+1,r,\ldots,r+k-1$) when we pick the
$a_{p_1},\ldots,a_{p_{2k+1}}$; such a choice then yields
$\{\bar{m}_n|n\leq 2k+1\}=\{1,1,\ldots,k,k,k+1\}$ and
\[ \sum_{n=1}^{2k+1}
\bar{m}_n=\frac{k(k+1)}{2}+\frac{(k+1)(k+2)}{2}=(k+1)^2.\]  If
$x_0\neq 0$, we have some $n\in [k+2,2k+1]$ with $p_n=0$ and so
$\bar{m}_n=n >k+1$; in this vein, one may easily check that \[
\sum_{n=1}^{2k+1} \bar{m}_n\geq (k+1)^2+\frac{x_0(x_0+1)}{2};\] in
particular \[ \sum_{n=1}^{2k+1} \bar{m}_n\geq (k+1)^2 +x_0,\] with
equality requiring that either $x_0=0$ and $\{\bar{m}_n|n\leq
2k+1\}=\{1,1,\ldots,k,k,k+1\}$ or $x_0=1$ and $\{\bar{m}_n|n\leq
2k+1\}=\{1,1,\ldots,k,k,k+2\}$.

Thus, \begin{align} (k+1)^2+x_0&\leq
\sum_{n=1}^{2k+1}\bar{m}_n=\sum_{n=1}^{2k+1}(n-p_n) \nonumber \\
&= (k+1)(2k+1)-\sum_{n=1}^{k+1}p_n-\sum_{n=k+2}^{2k+1}p_n
\nonumber \\
&=(k+1)^2-\sum_{j=1}^{k}jw_j+\sum_{n=k+2}^{2k+1}(k+1-p_n)
\nonumber \\ &\leq (k+1)^2-\sum_{j=1}^{k}jw_j+\sum_{n=k+2,p_n\leq
k}^{2k+1}(k+1-p_n) \nonumber \\
&=(k+1)^2-\sum_{j=1}^{k}jw_j+\sum_{j=0}^{k}(k+1-j)x_j \end{align}

So \begin{equation} kz_0+\sum_{j=1}^k(k+1-j)z_j\geq
kx_0+\sum_{j=1}^{k}(k+1-j)x_j\geq \sum_{j=1}^{k}jw_j\geq
\sum_{j=1}^{k}(j-1)w_j,\label{zx}\end{equation} \emph{i.e.},
$k\sum_{j=0}^{k}z_j+\sum_{j=1}^{k}(1-j)(w_j+z_j)\geq 0$, so that
since $\sum_{j=0}^{k}(w_j-1)=0$ and $i_j=w_j+z_j$, we at last
conclude that \begin{equation}\label{geq}
k\sum_{j=0}^{k}(i_j-1)+\sum_{j=1}^{k}(1-j)i_j \geq
0.\end{equation}  In light of Equation \ref{rearrange}, this shows
that $e(i_0,\ldots,i_r)\geq r+l^2-1$ with equality if and only if
equality holds in (\ref{geq}); equality in (\ref{geq}) requires
among other things that
\begin{itemize} \item[(i)] either $x_0=0$ and $\{\bar{m}_n|n\leq
2k+1\}=\{1,1,\ldots,k,k,k+1\}$ or $x_1=1$ and $\{\bar{m}_n|n\leq
2k+1\}=\{1,1,\ldots,k,k,k+2\}$; and \item[(ii)] due to (\ref{zx}),
$z_j=x_j$ for $j\leq k$ (so that for $j\leq k$ all of the $a_j$ in
our term $\prod_{j=0}^{r}a_{j}^{i_j}$ come from the first $2k+1$
columns of the resultant matrix). \end{itemize}

For $n=1,2,3,4$ let $M_n$ denote the $(2k+1)\times (2k+1)$ matrix
constructed from the resultant matrix $(\ref{res})$ by taking
columns $1$ through $2k+1$ and rows $1,\ldots,k,r,\ldots,r+k$ (for
$n=1$), rows $1,\ldots,k+1,r,\ldots,r+k-1$ (for $n=2$), rows
$1,\ldots,k,r,r+k-1,\ldots,r+k+1$ (for $n=3$), or rows
$1,\ldots,k,k+2,r,\ldots,r+k-1$ (for $n=4$).  Let $M'_n$ be the
$(2r-2k-2)\times (2r-2k-2)$ constructed from the other rows and
columns. Assume that our term $\prod_{j=0}^{r}a_{j}^{i_j}$ in the
resultant gives rise to the lowest possible value of
$e(i_0,\ldots,i_j)$. (i) above then ensures that
$\prod_{j=0}^{r}a_{j}^{i_j}$ is constructed by multiplying a term
in the determinant of one of the $M_n$ by a term in the
determinant of the corresponding $M'_n$.  In searching for the
optimal such monomial, we may then vary the contributions from
$M_n$ and $M'_n$ separately.  But on examining the form of the
$M_n$, one sees immediately that the term in $\det(M_n)$ giving
rise to the \emph{strictly} lowest possible power of $\ep$ is
obtained by a product of $k+1$ $a_0$'s (from columns 1 through
$k+1$ for $n=1,2$ and columns $1,\ldots,k,k+2$ for $n=3,4$) and
$k$ $a_{k+1}$'s (and in particular contains no $a_j$ for $1\leq
j\leq k$).  By (ii), any optimal monomial from $M'_n$ can't
contain any $a_j$ with $j\leq k$.  Thus any
$\prod_{j=1}^{r}a_{j}^{i_j}$ with $(i_1,\ldots,i_k)\neq
(0,\ldots,0)$ must have $e(i_1,\ldots,i_r)$ strictly greater than
the lowest possible value (which has been shown above to be
$n+l^2-1$).  This proves the lemma.\end{proof}

We now recall the linear algebra definition of the relative
Hilbert scheme from \cite{Smith}.  Let
\begin{equation}\label{hrmodel} \tilde{\mathcal{H}}_r=\{(A,B,t,v)\in
M_r(\mathbb{C})^2\times D^2\times\mathbb{C}^r| AB=BA=tId,
(*)\},\end{equation} where the stability condition (*) states that
the matrices A and B share no proper invariant subspaces
containing the vector $v$. The relative Hilbert scheme of the map
$(z,w)\mapsto zw$ is then
\[ \mathcal{H}_r=\tilde{\mathcal{H}}_r/GL_r(\mathbb{C}),\] where
$GL_{r}(\mathbb{C})$ acts by \[ g\cdot
(A,B,t,v)=(gAg^{-1},gBg^{-1},t,gv).\]  The projection map $F\co
\mathcal{H}_r\to D^2$ is just $[A,B,t,v]\mapsto t$. To briefly
motivate this, remark that a point of the $r$-fold relative
Hilbert scheme of $f$ is naturally viewed from an
algebro-geometric standpoint as an ideal $I\leq \mathbb{C}[z,w]$
with the property that $V=\mathbb{C}[z,w]/I$ is an $r$-dimensional
vector space and, for some $t$, $I$ is supported on $f^{-1}(t)$
(\emph{i.e.}, $\langle zw-t\rangle <I$).  To go from such an ideal
to an element of $\mathcal{H}_r$, let $v\in V$ be the image of
$1\in \mathbb{C}[z,w]$ under the projection, and let $A$ and $B$
be the operators on $V$ defined by multiplication by the
polynomials $z$ and $w$ respectively.  For more details see
\cite{N} and \cite{Smith}.

Given $[A,B,t,v]\in \mathcal{H}_r$, the fact that $A$ and $B$
commute implies that they can be simultaneously conjugated to be
upper triangular; assuming that this has been done, the natural
map $\phi_t\co F^{-1}(t)\to Sym^r f^{-1}(t)$ takes $[A,B,t,v]$ to
$\{(A_{11},B_{11}),\ldots,(A_{rr},B_{rr})\}$.  For $t\neq 0$,
according to (\ref{hrmodel}), $A$ is invertible and $B=tA^{-1}$, so
$\phi_t$ is an isomorphism; $\phi_0$, meanwhile, is a nontrivial
partial resolution.  On the diagonal $\Delta\subset
\mathcal{H}_r$, $A$ and $B$ will both have repeated eigenvalues,
occurring in corresponding Jordan blocks.

The main result of this section is:
\begin{thm} \label{appmain} Let $F\co \mathcal{H}_r\to D^2$ denote the $r$-fold
relative Hilbert scheme of the map $(z,w)\mapsto zw$, $\phi_0$ the
partial resolution map $F^{-1}(0)\to Sym^r\{zw=0\}$, and
$\Delta\subset \mathcal{H}_r$ the diagonal stratum. At any point
$p\in \Delta\cap F^{-1}(0)$ with
$\phi_0(p)=\{(0,0),\ldots,(0,0)\}$, where $T_p\Delta$ is the
tangent cone to $\Delta$ at $p$, we have $T_p\Delta\subset T_p
F^{-1}(0)$.
\end{thm}

\begin{proof}  According to the above description, the points $p$
under concern are of the form $[A,B,0,v]$ with $A$ and $B$ both
nilpotent matrices such that $AB=BA=0$  Further, letting $k$ be
such that $A^kv\neq 0$ but $A^{k+1}v=0$, the stability condition
$(*)$ in (\ref{hrmodel}) ensures that, where $r=k+l+1$, \[
\{A^kv,\ldots,Av,B^lv,\ldots,Bv,v\}\] is a basis for
$\mathbb{C}^r$.  All operators on $V\cong \mathbb{C}^r$ appearing
in the rest of the proof will be written as matrices in terms of
this basis.

Since $AB=0$ we can write \[
B^{l+1}v=aA^kv+\sum_{i=1}^{l}b_{l-i}B^iv.\]  With respect to our
above basis, we have \[ A=\left( \begin{array}{c c c c| c c c c
|c}0
& 1 & & & 0 & \cdot & \cdot & 0 & 0 \\ & \ddots & \ddots  & & \cdot & \cdot & &\cdot & \vdots \\  & & 0 & 1 & \cdot & & \cdot &\cdot & 0 \\
& & & 0 &  0 &\cdot & \cdot & 0 & 1 \\ \hline 0 & & \cdots & 0&0 &
\cdots & & & 0\\\vdots & & & & & & & &\vdots \\ 0 & & \cdots & 0&0
& \cdots & & & 0 \end{array}\right),\]
 \[ B=\left( \begin{array}{c c c c|c| c c
 c c}0&\cdot&\cdot&0&a&0&\cdot&\cdot&0\\
 \cdot &\cdot & &\cdot & 0&\cdot & & & \cdot\\
 \cdot& &\cdot&\cdot& \vdots&\cdot & &\ddots &\cdot \\
 \cdot & & & \cdot & 0 & 0 & \cdot &\cdot & 0\\
\hline 0 &\cdot &\cdot & 0 & b_1 & 1 & 0 &\cdots & 0\\
\hline \cdot &\cdot & & \cdot & \vdots & 0 & 1 &\cdots & 0\\
\cdot & & & \cdot & b_{l-1} & \cdot & \ddots & \ddots & \vdots \\
 \cdot& & \cdot & \cdot &b_l & \cdot & &\cdot & 1 \\
 0 & \cdot & \cdot & 0 & 0 & 0 &\cdot & \cdot & 0
 \end{array} \right), \] and $v=e_r=(0,\ldots,0,1)$ (in both of the above matrices, the upper left block is of size $k\times k$).
 Let \[
 (C,D,\mu,w)\in T_{(A,B,0,e_r)}\tilde{\mathcal{H}_r}.\]  Letting $\pi\co \tilde{\mathcal{H}}_r\to \mathcal{H}_r$
 be the projection, we have $\mu=F_*(\pi_*)_{(A,B,0,e_r)}(C,D,\mu,w)$, so our goal is to show that if $(C,D,\mu,w)$
 is tangent to $\pi^{-1}\Delta$ then $\mu=0$.  Linearizing the defining equations for $\tilde{H}_r$ gives
 \[CB+AD=BC+DA=\mu Id,\]
 which implies, among other things,
 \[ \mbox{For $i>k$, }\begin{cases} aC_{i1}+\sum b_m
 C_{i,k+m}=\mu\delta_{i,j} \\ C_{i,j-1}=\mu\delta_{i,j} \mbox{ if
 } j\geq k+2\end{cases} \]\[\mbox{For $j=1$ or $k+1\leq j\leq
 r-1$, } \begin{cases} aC_{k+1,j}=\mu\delta_{1,j} \\
 b_{i-k}C_{k+1,j}+C_{i+1,j}=\mu\delta_{i,j} \mbox{ if } k+1\leq
 i\leq r-1\end{cases} \]

 If $a=0$, we have $\mu=\mu\delta_{1,1}=aC_{k+1,1}=0$ and we are
 done.

If $a\neq 0$, we find from the above equations that\[ C=\left(
\begin{array}{c c c c| c c c c |c}
* & * &* &* & *& * & * & * & * \\
*& *& *& * & * & * &* &*& * \\
* & *&* & * & * &* & *&*& *\\
* &* & *& * & * &*& *& * & * \\\hline
\mu/a & * & *& *&0 & 0 &\cdot &0 & *\\
\hline -b_1\mu/a & * & *& *&\mu &0 &\cdots & 0&* \\
-b_2\mu/a &* &* & *&0 & \mu& \cdot & 0 & * \\
\vdots & * & * & * & \cdot & \cdot & \ddots& \cdot & * \\
-b_l\mu/a & * & * & * & 0 & \cdot & \cdot &\mu & *
\end{array}\right), \] where again the upper left block is size
$k\times k$ and all asterisks denote undetermined entries.

We consider now the characteristic polynomials of the matrices
$A+\ep C$ for small $\ep$. The matrix $A+\ep C-\lambda Id$ is
\[\left(
\begin{array}{c c c c| c c c c |c}
-\lambda+\ep C_{11} & 1+\ep C_{12} &* &* & *& * & * & * & * \\
*& -\lambda+\ep C_{22}& \ddots& * & * & * &* &*& * \\
* & *&\ddots & 1+\ep C_{k-1,k} & * &* & *&*& *\\
* &* & *& -\lambda+\ep C_{kk} & * &*& *& * & 1+\ep C_{rk} \\\hline
\ep\mu/a & * & *& *&-\lambda & 0 &\cdot &0 & *\\
\hline -\ep b_1\mu/a & * & *& *&\ep\mu &-\lambda &\ddots & 0&* \\
-\ep b_2\mu/a &* &* & *&0 & \ep\mu& \cdot & 0 & * \\
\vdots & * & * & * & \vdots & \cdot & \ddots& -\lambda & * \\
-\ep b_l\mu/a & * & * & * & 0 & \cdot & \cdot &\ep\mu &
-\lambda+\ep C_{rr} \end{array}\right)\] where an asterisk in the
$(i,j)$th entry signifies $\ep C_{ij}$. When we expand the
determinant of this matrix, among the terms that we obtain are \[
(-\lambda)^r \mbox{  and  } \ \pm (-\ep b_m\mu/a)\cdot
1^{k-1}(-\lambda)^{m-1}(\ep\mu)^{l-m+1}=\pm
\frac{\mu^{l-m+2}b_m}{a}\ep^{l-m+2}\lambda^{m-1}; \] note that
these latter have combined degree exactly $l+1$ in $\ep$ and
$\lambda$. Any other term in the expansion of the determinant will
have degree at least 1 in $\ep$ and at least $l+2$ in $\ep$ and
$\lambda$ combined, the reason being that each of the entries
denoted with an asterisk above lies in either the same row or the
same column as an entry of form $1+\ep C_{ij}$, so a term in the
determinant containing one of the asterisked entries can contain
at most $k-1$ of the $k$ $(1+\ep C_{ij})$'s and hence must contain
at least $r-(k-1)=l+2$ other terms, each of which is of combined
order at least 1 in $\ep$ and $\lambda$.  In other words, for
constants $c_1,\ldots,c_{k+l+1}$ where
\[ c_{k+m}=\pm \frac{b_{l+1-m}\mu^m}{a} \mbox{ for $1\leq m\leq l$
and } c_{k+l+1}=\frac{\mu^{l+1}}{a},\] the characteristic
polynomial of $A+\ep C$ has form \[ p_{A+\ep
C}(x)=(-x)^r+\ep\sum_{a=1}^{k}(c_a+O(\ep))(-x)^{r-a}+\sum_{b=1}^{l+1}\ep^b(c_{k+b}+O(\ep))(-x)^{l+1-b},\]
which, since $r=k+l+1$, is precisely the sort of polynomial
considered in Lemma \ref{discrim}.  By replacing $\ep$ with
$\nu\ep$ in the statement of that lemma, we see that the
polynomial $P(a_{r-l},\ldots,a_r)$ provided by its conclusion
scales as \[ P(\nu a_{r-l},\nu^2
a_{r-l+1},\ldots,\nu^{l+1}a_r)=\nu^{r+l^2-1}P(a_{r-l},\ldots,a_r),\]
so that \begin{align*}
P(c_{k+1},\ldots,c_r)&=P\left(\pm\frac{b_l\mu}{a},\pm\frac{b_{l-1}\mu^2}{a},\ldots,\pm\frac{b_1\mu^l}{a},\frac{\mu^{l+1}}{a}\right)\\
&=\mu^{r+l^2-1}P(\pm b_l/a,\pm b_{l-1}/a,\ldots,\pm
b_1/a,1/a).\end{align*}  So since $P$ is not the zero polynomial,
at least for a generic initial choice of our base point
$[A,B,0,e_r]$ (equivalently, for generic $a,b_1,\ldots,b_l$) we
conclude that if $(C,D,\mu,w)\in
T_{(A,B,0,e_r)}\tilde{\mathcal{H}}_r$, we have \begin{equation}
\label{conc} \delta(p_{A+\ep
C})=\mu^{r+l^2-1}M\ep^{n+l^2-1}+O(\ep^{n+l^2}),\end{equation}
where $M$ is a nonzero constant depending only on $A$.  Let
\begin{align*} \tilde{\Delta}_1&=\{(A',B',t,v')\in
\tilde{\mathcal{H}}_r|A' \mbox{ has a repeated eigenvalue}\}\\
&=\{(A',B',t,v')\in
\tilde{\mathcal{H}}_r|\delta(p_A)=0\}\end{align*} Equation
\ref{conc} then shows that, for $(C,D,\mu,w)\in
T_{(A,B,0,e_r)}\tilde{\mathcal{H}}_r$, \[ (C,D,\mu,w)\in
T_{(A,B,0,e_r)}\tilde{\Delta}_1\Leftrightarrow \mu=0\]
($T_{(A,B,0,e_r)}\tilde{\Delta}_1$ denoting the tangent cone at
$(A,B,0,e_r)$ if $\tilde{\Delta}_1$ is singular there).  Where
again $\pi\co \tilde{\mathcal{H}}_r\to\mathcal{H}_r$ is the
projection, we have $T\Delta\subset \pi_*T\tilde{\Delta}_1$, so if
$\alpha\in T_{[A,B,0,e_r]}\Delta$, writing
$\alpha=\pi_*(C,D,\mu,w)$, we have that $F_*\alpha=\mu=0$.  This
conclusion initially only applies at those $[A,B,0,e_r]\in \Delta$
which are generic in the sense that $P(\pm b_l/a,\pm
b_{l-1}/a,\ldots,\pm b_1/a,1/a)\neq 0$, but then since the
conclusion is a closed condition it in fact applies to all
$[A,B,0,e_r]$ lying on the diagonal $\Delta$. \end{proof}

\end{document}